\documentclass[brochure, 12pt,french]{bourbaki}

\usepackage{pgfpict2e}
\usepackage{amssymb,amsfonts,amsmath,footnote}
\usepackage{epsfig}
\usepackage{graphics}
\usepackage{tikz}
\usepackage{mathrsfs}

\usepackage{indentfirst}

\usepackage[hmargin=20mm,vmargin=20mm]{geometry}
\usepackage[applemac]{inputenc}
\usepackage[T1]{fontenc}
\usepackage{lmodern}

\usepackage[francais]{babel}
\newtheorem{theorem}{\indent Th\'{e}or\`{e}me}
\newtheorem{rem}[theorem]{\indent Remarque}
\newtheorem{lem}[theorem]{\indent Lemme}

\newtheorem{cor}[theorem]{\indent Corollaire}

 \newcommand{\R}{{\mathbb R}}
\newcommand{\C}{{\mathbb C}}

\newcommand{\Z}{\mathbb Z}

\newcommand{\Hyp}{{\mathbb H}}

\newcommand{\AAA}{{\mathcal A}}
\newcommand{\LLL}{{\mathcal L}}
\newcommand{\HHH}{{\mathcal H}}
\newcommand{\CCC}{{\mathcal C}}

\date{Octobre 2012}
\bbkannee{65\`eme ann\'ee, 2012-2013}
\bbknumero{1060}
\title{Les exposants de Liapounoff du flot de Teichm\"uller}
\subtitle{d'apr\`es Eskin-Kontsevich-Zorich}
\author{Julien GRIVAUX \& Pascal HUBERT  }

\address{Universit\'{e} d'Aix-Marseille\\
Centre de Math\'ematiques et Informatique (CMI)\\
Technop\^ole Ch\^ateau-Gombert\\
39, rue F. Joliot Curie\\
F--13453 Marseille Cedex 13}
\email{jgrivaux@cmi.univ-mrs.fr}
\email{hubert@cmi.univ-mrs.fr}

\begin{document}

\maketitle

Dans ce texte, nous pr\'{e}sentons un r\'{e}sultat d\^{u} \`{a} Eskin, Kontsevich et Zorich \cite{EKZ2} concernant les exposants de Liapounoff du flot de Teichm\"{u}ller. 

\section{Introduction}

\subsection{Surfaces de translation et 1--formes holomorphes}
Nous allons en premier lieu fixer le cadre g\'{e}n\'{e}ral (pour des textes introductifs, on renvoie le lecteur aux r\'{e}f\'{e}rences suivantes: \cite{MT}, \cite{Vi}, \cite{Yo1}, \cite{Yo2}, \cite{Zo4}).
Une surface de translation compacte est la donn\'{e}e d'une 1--forme holomorphe globale non nulle $\omega$ sur une surface de Riemann compacte $X$ \footnote{Eskin-Kontsevich-Zorich traitent aussi le cas des formes diff\'{e}rentielles quadratiques qui est analogue, mais nous nous limiterons aux formes holomorphes pour simplifier l'exposition.}.
G\'{e}om\'{e}triquement, une telle surface se repr\'{e}sente comme un polygone dans le plan complexe dont on a identifi\'{e} par translation des c\^{o}t\'{e}s parall\`{e}les et de m\^{e}me longueur. La 1--forme $\omega$ sur $X$ est induite par $dz$, et ses z\'{e}ros sont certains sommets du polygone. De plus, $X$ est munie d'une m\'{e}trique plate h\'{e}rit\'{e}e de la m\'{e}trique plate naturelle du polygone, avec des singularit\'{e}s coniques sur le lieu d'annulation de $\omega$. De mani\`{e}re pr\'{e}cise, un z\'{e}ro d'ordre $k$ 
correspond \`{a} un point conique d'angle $2(k+1)\pi$. 
\par \medskip
L'aire de $X$ pour la m\'{e}trique plate est l'aire euclidienne d'un polygone associ\'{e}, et on a
\[
\mathrm{Aire}(X)=\displaystyle{\frac{i}{2}\int_X \omega \wedge \overline{\omega}}. 
 \]
Enfin le genre $g$ de $X$ est d\'{e}termin\'{e} par la donn\'{e}e combinatoire $(k_1, \dots, k_r)$ de $(X, \omega)$, qui est la liste des ordres de multiplicit\'{e}s des z\'{e}ros de $\omega$, via la formule
\[
\sum_{i=1}^r k_i = 2g-2.
\]
L'exemple de base de surface de translation est le tore plat qui est \'{e}videmment une surface de translation sans singularit\'{e}\footnote{Pour \^{e}tre coh\'{e}rent avec le reste de la th\'{e}orie, l'origine du tore est tout de m\^{e}me un point marqu\'{e}.}, le polygone associ\'{e} \'{e}tant un parall\'{e}logramme dont les c\^ot\'es oppos\'{e}s sont identifi\'{e}s.

 \begin{figure}[ht]
 \begin{center}
 \begin{tikzpicture}[scale=3]
 \draw[very thick]  (0,0)--(1,0);
 \draw[very thick]  (1,0)--(1,1);
 \draw[very thick]  (1,1)--(0,1);
 \draw[very thick]  (0,1)--(0,0);
 \draw(0.25,0)--(1, 0.325);
 \draw(0, 0.325,0)--(1, 0.73);
 \draw(0, 0.73)--(0.65,1);
 \end{tikzpicture}
 \end{center}
  \caption{Tore plat et flot lin\'{e}aire. }
 \label{fig:tore}
  
  \end{figure}
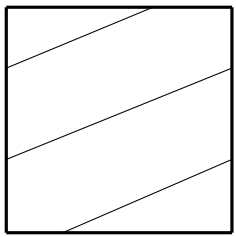

Pour  tout angle $\theta$, le flot lin\'{e}aire de direction $\theta$ est bien d\'{e}fini sur $X$, et ses propri\'{e}t\'{e}s dynamiques sont d'un grand int\'{e}r\^{e}t. 

\subsection{Espace des modules des 1--formes holomorphes}

Un genre $g$ \'{e}tant fix\'{e}, l'espace des modules des 1--formes holomorphes est l'ensemble des surfaces de translation compactes $(X, \omega)$ de genre $g$ modulo l'action naturelle des diff\'{e}omorphismes. Cet espace de modules est naturellement stratifi\'{e} par les donn\'{e}es combinatoires. Pour toute donn\'{e}e combinatoire $(k_1, \dots, k_r)$ nous noterons $\HHH(k_1, \dots, k_r)$ la strate correspondante. On peut montrer que $\HHH(k_1, \dots, k_r)$
est une orbifolde complexe de dimension $2g+r-1$ localement model\'{e}e sur le groupe de cohomologie relative ${\mathrm{H}}^1(X, \Sigma, \C)$, o\`{u} $\Sigma$ est le lieu d'annulation de $\omega$. 
\par \medskip
On rappelle bri\`{e}vement la construction de coordonn\'{e}es locales orbifoldes sur $\HHH(k_1, \dots, k_r)$.
Fixons une base symplectique $(A_i, B_i)_{i = 1, \dots, g}$ de ${\mathrm{H}}_1(X,\Z)$,  et soient $(C_i)_{i= 1, \dots, r-1}$ des chemins reliant un z\'{e}ro fix\'{e}  de $\omega$ \`{a} tous les autres.
Les coordonn\'{e}es locales sont les p\'{e}riodes de $\omega$ le long de ces chemins, c'est-\`{a}-dire les int\'{e}grales
\[
\int_{A_1} \omega, \dots, \int_{A_g} \omega, \int_{B_1} \omega, \dots, \int_{B_g} \omega,  \int_{C_1} \omega, \dots, \int_{C_{r-1}} \omega.
\]
Les changements de cartes correspondent \`{a} un changement de base symplectique et sont donc lin\'{e}aires, ce qui entra\^{\i}ne que chaque strate $\HHH(k_1, \dots, k_r)$ admet une structure affine. De mani\`{e}re plus pr\'{e}cise, 
dans les coordonn\'{e}es ci-dessus, une matrice de changement de base est de la forme 
\[
\begin{pmatrix}
U & V \\
0 & I_{r-1} 
\end{pmatrix}
\]
o\`{u} $U$ est dans $\textrm{Sp}(2g, \Z)$ et $V$ est dans $\textrm{M}_{2g \times (r-1)} (\Z )$. Une telle matrice est \'{e}videmment dans $\textrm{SL}_{2g+r-1}(\Z)$. Par cons\'{e}quent,
la mesure de Lebesgue de $ \C^{2g+r-1}$ est bien d\'{e}finie globalement sur la strate $\HHH(k_1, \dots, k_r)$. Le r\'eseau entier $(\Z + i\Z)^{2g+r-1} = H^1(X,\Sigma, \Z+i\Z)$ est invariant par changement de cartes ; ce r\'eseau nous fournit une normalisation {\it naturelle} pour la mesure de Lebesgue: on demande que son covolume soit $1$. 
\par \medskip
Le groupe $\textrm{GL}_2(\R)$ agit de fa\c{c}on naturelle sur $\HHH(k_1, \dots, k_r).$ Cette action correspond \`{a} l'action lin\'{e}aire sur les polygones, qui est bien d\'{e}finie au niveau des surfaces de translation et pr\'{e}serve la donn\'{e}e combinatoire. Dans les coordonn\'{e}es des p\'{e}riodes introduites dans la section pr\'{e}c\'{e}dente, l'action de toute matrice $M$ de $\textrm{GL}_2(\R)$ est diagonale~: 
\begin{align*}
&M \cdot \Bigl(\int_{A_1} \omega, \dots, \int_{A_g} \omega, \int_{B_1} \omega, \dots, \int_{B_g} \omega,  \int_{C_1} \omega, \dots, \int_{C_{r-1}} \omega \Bigr) \\
&\qquad = \Bigl(M \cdot \int_{A_1} \omega, \dots, M \cdot \int_{A_g} \omega, M \cdot \int_{B_1} \omega, \dots, M \cdot \int_{B_g} \omega,  M \cdot \int_{C_1} \omega, \dots, M \cdot\int_{C_{r-1}} \omega\Bigr).
\end{align*}

\noindent Cette action est  $\R$--lin\'{e}aire dans chaque coordonn\'{e}e complexe (vue comme coordonn\'{e}e \`{a} valeurs dans $\R^2$).
\par \medskip
L'action de $\textrm{GL}_2(\R)$ permet de feuilleter les strates par les orbites, ce feuilletage sera particuli\`{e}rement important dans la suite.
Il est \'{e}vident par la formule pr\'{e}c\'{e}dente que le groupe $\textrm{SL}_2(\R)$  pr\'{e}serve la mesure de Lebesgue de $\HHH{(k_1, \dots, k_r)}$, ainsi que l'aire des surfaces de translation.
Notons $\HHH_1(k_1, \dots, k_r)$ la sous-vari\'{e}t\'{e} de codimension r\'{e}elle un form\'{e}e des surfaces de translation d'aire 1 et de donn\'{e}e combinatoire $(k_1, \dots, k_r)$.
La mesure de Lebesgue sur $\HHH{(k_1, \dots, k_r)}$ induit une mesure $\nu$ sur $\HHH_1{(k_1, \dots, k_r)}$. Dans les coordonn\'{e}es des p\'{e}riodes, la mesure d'un ensemble de $\HHH_1{(k_1, \dots, k_r)}$ est celle du c\^{o}ne de sommet l'origine bord\'{e} par cet ensemble à un facteur dimensionnel pr\`es.
\par \medskip
Un sous-groupe \`{a} un param\`{e}tre de $\textrm{SL}_2(\R)$ joue un r\^{o}le fondamental dans toute la suite, c'est le groupe des matrices diagonales
\[
\begin{pmatrix}
e^t & 0 \\
0 & e^{-t}
\end{pmatrix}
\]
\par \smallskip
Le flot associ\'{e} est le flot de Teichm\"{u}ller, qui est le flot g\'{e}od\'{e}sique pour la m\'{e}trique de Teichm\"{u}ller.
Un des premiers th\'{e}or\`{e}mes importants concernant ce flot est d\^{u} ind\'{e}pendamment \`{a} Masur et Veech en 1982:

\begin{theorem}[\cite{Ma1}, \cite{Ve1}]
La mesure de Lebesgue $\nu$ sur $\HHH_1{(k_1, \dots, k_r)}$ est finie. De plus, le flot g\'{e}od\'{e}sique de Teichm\"{u}ller est ergodique sur chaque composante connexe de 
$\HHH_1{(k_1, \dots, k_r)}$.
\end{theorem}
Les composantes connexes de $\HHH_1{(k_1, \dots, k_r)}$ ont \'{e}t\'{e} classifi\'{e}es par Kontsevich et Zorich dans \cite{KZ2}, il y en a au plus 3. Dans le cas des diff\'erentielles quadratiques, elles ont \'et\'e classifi\'ees par Lanneau \cite{La}.\begin{rem}
Dans toute la suite, nous parlerons de strate l\`a o\`u il faudrait parler en toute rigueur de composante connexe de strate, ceci pour ne pas alourdir le texte. La mesure $\nu$ d\'epend de la composante connexe consid\'er\'ee~: c'est son support. 
\end{rem}
 
\subsection{Fibr\'{e} de Hodge et Cocycle de Kontsevich-Zorich} \label{subsection:Hodge}

\'{E}tant donn\'{e} une surface de Riemann $X$ de genre $g$, rappelons la d\'{e}finition de la norme de Hodge sur l'espace ${\mathrm{H}}^1(X,\R)$. 
La d\'{e}composition de Hodge s'\'{e}crit \[
{\mathrm{H}}^1(X,\C) = {\mathrm{H}}^{1,0}(X) \oplus {\mathrm{H}}^{0,1}(X)
\]
 o\`{u} ${\mathrm{H}}^{1,0}(X)$ est l'espace vectoriel de dimension $g$ des 1--formes holomorphes sur $X$ et ${\mathrm{H}}^{0,1}(X)$ l'espace vectoriel  des 1--formes anti-holomorphes. 
La forme d'intersection sur ${\mathrm{H}}^1(X,\C)$ donn\'{e}e par 
\[
\iota(\omega_1,\omega_2) = \frac{i}{2} \int_X \omega_1 \wedge \overline{\omega}_2
\]
est d\'{e}finie positive sur le sous-espace ${\mathrm{H}}^{1,0}(X)$ et d\'{e}finie n\'{e}gative sur le sous-espace conjugu\'{e} ${\mathrm{H}}^{0,1}(X)$. Notons 
$\phi \colon {\mathrm{H}}^{1,0}(X) \rightarrow {\mathrm{H}}^{1}(X,\R)$ l'isomorphisme qui \`{a} toute forme holomorphe associe la classe de sa partie r\'{e}elle.
C'est un isomorphisme de $\R$--espaces vectoriels qui permet de munir ${\mathrm{H}}^{1}(X,\R)$ de la norme de Hodge~:
si $\alpha$ est dans ${\mathrm{H}}^{1}(X,\R)$, on pose
\[
\vert \vert \alpha \vert \vert^2_{\textrm{\tiny{Hodge}}} = \iota \bigl(\phi^{-1} (\alpha), \phi^{-1} (\alpha)\bigr).
\]

On consid\`{e}re maintenant le fibr\'{e} de Hodge r\'{e}el au-dessus d'une strate $\HHH{(k_1, \dots, k_r)}$ dont la fibre au dessus de chaque point $(X, \omega)$ est ${\mathrm{H}}^1(X,\R)$. Ce fibr\'{e} de Hodge provient d'un syst\`{e}me local de $\R$-espaces vectoriels, ce qui permet
d'identifier localement les fibres voisines entre elles. L'identification se fait gr\^{a}ce \`{a} la connexion plate sur le fibr\'{e} de Hodge associ\'{e}e \`{a} ce syst\`{e}me local, appel\'{e}e connexion de Gau\ss-Manin. Une grande partie de la g\'{e}om\'{e}trie du probl\`{e}me est contenue dans le fait que le fibr\'{e} de Hodge complexe de fibre ${\mathrm{H}}^1(X,\C)$ est un fibr\'{e} plat, mais que le sous-fibr\'{e} de fibre 
${\mathrm{H}}^{1,0}(X)$ qui tient compte de la structure complexe de $X$ ne respecte pas cette structure plate (c.a.d. n'est pas invariant par la connexion de Gau\ss-Manin).
\par \medskip
Etant donn\'{e} une surface de translation $X$ et une classe de cohomologie $\alpha$ dans $\mathrm{H}^1(X, \R)$, on souhaite comprendre la croissance de la norme de Hodge $||\alpha_t||_{g_t X}$ quand $t$ tend vers l'infini, o\`{u} $g_t$ est le flot de Teichm\"{u}ller et $\alpha_t$ est la classe transport\'{e}e parall\`{e}lement le long du flot \`{a} partir de $\alpha$. 
\par \medskip
La monodromie de la connexion de Gau\ss-Manin d\'{e}finit un cocycle: en tout point $X$ de la strate, on dispose d'une repr\'{e}sentation du groupe fondamental de la composante connexe de $X$ dans $\mathrm{H}^1(X, \R)$ qui est bien d\'{e}finie \`{a} conjugaison pr\`{e}s. Par le th\'{e}or\`{e}me d'Ehresmann, l'action d'un \'{e}l\'{e}ment du groupe fondamental est celle d'un diff\'{e}omorphisme orient\'{e} sur l'homologie; c'est donc une matrice symplectique car tout diff\'{e}omorphisme orient\'{e} pr\'{e}serve la forme d'intersection.
\par \medskip 
Expliquons intuitivement comment obtenir ce cocycle de fa\c{c}on concr\`{e}te le long des g\'{e}od\'{e}siques du flot de Teichm\"{u}ller. On fixe un petit ouvert $U$ de la strate dans lequel on peut identifier de mani\`{e}re canonique les espaces de cohomologie $\mathrm{H}^1(Y, \R)$ entre eux lorsque $Y$ parcourt $U$. Soit $X$ un point de $U$ et $\alpha$ une classe de cohomologie dans $\mathrm{H}^1(X, \R)$. On suit par transport parall\`{e}le la classe $\alpha$ sous l'action du flot de Teichm\"{u}ller $g_t$ jusqu'\`{a} revenir dans $U$ (ce qui se produit en temps fini pour presque tout $X$ car le flot est ergodique). 
La classe $\alpha_t$ obtenue dans $\mathrm{H}^1(g_t X, \R)$ s'identifie de mani\`{e}re canonique \`{a} une classe dans $\mathrm{H}^1(X, \R)$ qui est pr\'{e}cis\'{e}ment l'action du cocycle de Kontsevich-Zorich $G_t^{{\mathrm{KZ}}} $ sur la classe $\alpha$.
\subsection{Exposants de Liapounoff du flot de Teichm\"{u}ller} \label{exposants}
Comme le flot de Teichm\"{u}ller est ergodique, on peut appliquer le th\'{e}or\`{e}me d'Oseledets au cocycle de Kontsevich-Zorich. Il existe donc des nombres r\'{e}els $\lambda_1 > \dots > \lambda_k$  et  une d\'{e}composition 
\[
{\mathrm{H}}^1(X, \R) = E_1(\omega) \oplus \dots \oplus E_k(\omega)
\]
d\'{e}pendant mesurablement de $(X,\omega)$ telle que
si $\alpha$ est dans $E_i(\omega)$, on a 
\[
\lim_{t \to \pm \infty} \frac{1}{t} \log \, \vert \vert \, G_t^{{\mathrm{KZ}}} (\alpha)\, \vert \vert = \lambda_i.
\]

Il est plus ais\'{e} pour la suite de consid\'{e}rer  $2g$ exposants de Liapounoff $\lambda_1 \geq \lambda_2 \geq  \dots  \geq \lambda_{2g}$, chaque exposant \'{e}tant r\'{e}p\'{e}t\'{e} avec une multiplicit\'{e} correspondant \`{a} la dimension de l'espace d'Oseledets associ\'{e}. Comme le cocycle est symplectique, ces exposants v\'{e}rifient la relation de sym\'{e}trie~
$\lambda_{2g-i+1} = - \lambda_i$. De plus, on montre que $\lambda_1=1$ de la mani\`{e}re suivante: le sous-fibr\'{e} de rang $2$ du fibr\'{e} de Hodge dont la fibre en tout point $(X, \omega)$ est le plan r\'{e}el engendr\'{e} par $\mathrm{Re}\,\omega$ et  $\mathrm{Im}\,\omega$ est stable par le flot de Teichm\"{u}ller. Le cocycle de Kontsevich-Zorich restreint \`{a} ce plan n'est autre que l'action lin\'{e}aire donn\'{e}e par le flot g\'{e}od\'{e}sique de $\textrm{SL}_2 (\R)$, et les exposants de Liapounoff associ\'{e}s sont $1$ et $-1$. Ces exposants sont extr\'{e}maux d'apr\`{e}s le th\'{e}or\`{e}me de Teichm\"{u}ller.
\par \medskip
Les valeurs des autres exposants sont dans la plupart des cas totalement inconnues.
Forni \cite{Fo} a prouv\'{e} que $\lambda_g$ est strictement positif et Avila-Viana \cite{AV} que tous les exposants sont distincts. Ces r\'{e}sultats difficiles sont vrais uniquement dans le cas des strates et pour la mesure $\nu$. 
\par\medskip
\`{A} la suite d'exp\'{e}riences num\'{e}riques obtenues par l'algorithme de Rauzy-Veech (qui permet de discr\'{e}tiser le flot de Teichm\"{u}ller), Kontsevich et Zorich (\cite{Ko}, \cite{KZ1}) ont conjectur\'{e} vers le milieu des ann\'{e}es 90 que la somme des exposants positifs $\lambda_1 + \dots +\lambda_g$ est un nombre rationnel, ce qui est remarquable et a priori tr\`{e}s surprenant vu la d\'{e}finition des exposants. L'article d'Eskin-Kontsevich-Zorich donne une formule explicite pour cette somme. Combin\'{e} avec des r\'{e}sultats ant\'{e}rieurs, cette formule implique la rationalit\'{e} de la somme des exposants positifs.
\par\medskip
L'\'{e}tude de ces exposants de Liapounoff se justifie par le fait que le flot de Teichm\"{u}ller joue le r\^{o}le d'op\'{e}rateur de renormalisation pour les surfaces de translation, ce qui est connu depuis les travaux de Masur et Veech. Le th\'{e}or\`{e}me de Kerckhoff, Masur et Smillie \cite{KMS} affirme que pour toute surface de translation, le flot {\it lin\'{e}aire} est  uniquement ergodique dans presque toute direction. Supposons pour simplifier que le flot vertical est uniquement ergodique, fixons un long morceau d'orbite de ce flot que l'on ferme par un chemin de longueur born\'{e}e, et notons $\gamma_t$ cette courbe ferm\'{e}e. Le th\'{e}or\`{e}me ergodique nous assure que la quantit\'{e} $\,t^{-1} [\gamma_t]$ converge dans $\mathrm{H}_1(X, \R)$ \footnote{Cette limite s'appelle cycle asymptotique de Schwartzman.}. Zorich et  Forni (\cite{Zo2}, \cite{Zo3}, \cite{Fo}) ont montr\'{e} que les d\'{e}viations par rapport \`{a} cette moyenne sont gouvern\'{e}es par les exposants de Liapounoff d\'{e}crits ci-dessus. 
\par \medskip
Notons $(F_i(\omega))_{1 \leq i \leq 2g}$ la filtration d\'{e}croissante d'Oseledets pour le cocycle de Kontsevich-Zorich agissant sur $\mathrm{H}_1(X, \R)$. 
Zorich (\cite{Zo2}, \cite{Zo3})  montre que pour une surface de translation g\'{e}n\'{e}rique $(X, \omega)$ et tout \'{e}l\'{e}ment 
$f$ dans $F_i(\omega) \setminus F_{i+1}(\omega)$, on a 
\[
\limsup_{t \to +\infty} \frac{\log\vert \langle f, \gamma_t \rangle\vert}{\log t}  = \lambda_i
\]
lorsque $ \lambda_i >0$ (c'est-\`{a}-dire si $i \leq g$). De plus, la quantit\'{e}  $\vert \langle f, \gamma_t \rangle\vert$ est born\'{e}e si $f$ appartient \`{a} $F_{g+1}$.
Ce r\'{e}sultat est utilis\'{e} dans Delecroix-Hubert-Leli\`{e}vre \cite{DHL} pour comprendre la vitesse de diffusion d'un billard polygonal non compact, appel\'{e} mod\`{e}le windtree. Ici, non seulement l'existence des exposants mais aussi leur valeur est importante. 
\par\medskip
Pour terminer ces motivations, mentionnons que
le fait que $\lambda_g$ soit strictement positif pour les strates est un \'{e}l\'{e}ment essentiel de la preuve du th\'{e}or\`{e}me d'Avila et Forni \cite{AF} sur le m\'{e}lange faible des \'{e}changes d'intervalles.

\subsection{Constantes de Siegel-Veech}

Pour \'{e}noncer la formule d'Eskin-Kontsevich-Zorich sur la somme des exposants de Liapounoff du flot de Teichm\"{u}ller, il est n\'{e}cessaire d'introduire au pr\'{e}alable les constantes de Siegel-Veech.
\par \medskip
Rappelons que sur une surface de translation, on appelle lien de selles un segment g\'{e}od\'{e}sique (pour la m\'{e}trique plate) reliant deux singularit\'{e}s et n'en contenant aucune dans son int\'{e}rieur. Dans une surface de translation, les orbites p\'{e}riodiques arrivent par familles; elles forment des
cylindres bord\'{e}s par des liens de selles.  Pour un cylindre $\CCC$, on notera $w(\CCC)$ son p\'{e}rim\`{e}tre et $h(\CCC)$ sa hauteur. L'aire du cylindre est $w(\CCC) h(\CCC)$, son module $\mathrm{Mod}(\CCC)$ est $\dfrac{h(\CCC)^{\vphantom{\bigl(}}}{w(\CCC)}\,$.
 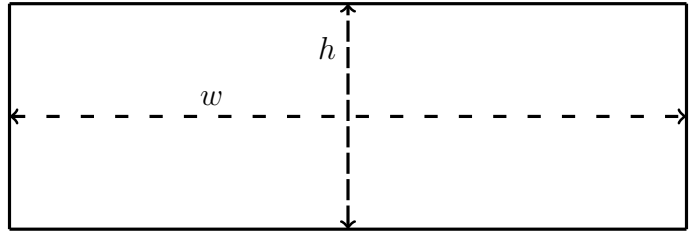
\begin{figure}[!h]
 \begin{center}
 \begin{tikzpicture}[scale=3]
  \tikzstyle{a line}=[very thick, dash pattern=on 8pt off 2pt];
   \tikzstyle{b line}=[very thick, dash pattern=on 5pt off 8pt];
 \draw[very thick]  (0,0)--(3,0);
 \draw[very thick]  (3,0)--(3,1);
 \draw[b line, arrows=<->]   (0,0.5)--(1,0.5) node [above left] {$w$}--(3,0.5);
 \draw[very thick]  (3,1)--(0,1);
 \draw[very thick]  (0,1)--(0,0);
   \draw[a line, arrows=<->] (1.5,0)--(1.5,0.7) node [above left] {$h$}--(1.5, 1);
 \end{tikzpicture}
 \end{center}
  \caption{Cylindre horizontal (les c\^{o}t\'{e}s verticaux sont identifi\'{e}s). }
 \label{fig:tore}
  \end{figure}
  \par
Pour toute surface de translation $(X, \omega)$ et tout r\'{e}el positif $T$, notons $N(X,T) $ le nombre de liens de selles de longueur inf\'{e}rieure ou \'{e}gale \`{a} $T$ sur $X$.
Masur \cite{Ma2} montre que cette quantit\'{e} a une croissance quadratique~: il existe une constante $C = C(X) >1$ telle que pour tout $T$ positif, 
\[
C^{-1}\, {T^2} \leq N(X,T) \leq C\, T^2.
\]
On peut se demander si $T^{-2} {N(X,T)}$ a une limite quand $T$ tend vers l'infini. Dans le cas du tore, il est bien connu que cette limite existe et vaut $6/\pi^2$ (qui est le terme dominant dans le probl\`{e}me du cercle). 
La question n'est toujours pas r\'{e}solue pour une surface de translation arbitraire. N\'{e}anmoins,
Eskin et Masur \cite{EM} montrent qu'il existe une constante $C(\nu)$ telle que pour $\nu$-presque toute surface de la strate $\HHH_1{(k_1, \dots, k_r)}$ on a
\[
\displaystyle\lim_{T \to +\infty} T^{-2}\, {N(X,T)}= \pi \,C(\nu).
\]
La constante $C(\nu)$ est appel\'{e}e constante de Siegel-Veech \footnote{Veech \cite{Ve3} a \'{e}t\'{e} le premier \`{a} travailler sur ce probl\`{e}me en adaptant aux surfaces de translation des id\'{e}es d\'{e}velopp\'{e}es par Siegel pour \'{e}tudier les espaces de r\'{e}seaux.}.   
\par \medskip
Pour ce qui suit, on va consid\'{e}rer un probl\`{e}me l\'{e}g\`{e}rement diff\'{e}rent: on s'int\'{e}resse au comptage de cylindres \textit{pond\'{e}r\'{e}s}, le poids attach\'{e} \`{a} chaque cylindre \'{e}tant son aire. Ceci peut \^{e}tre justifi\'{e} de mani\`{e}re g\'{e}om\'{e}trique par le fait que les cylindres sont des g\'{e}od\'{e}siques \og \'{e}paisses\fg{}, ils ont donc un poids. Pour tout r\'{e}el positif $R$, on pose 
\[
N_{\mathrm{aire}}(X, R) = \frac{1}{\mathrm{Aire(X)}} \,\sum_{\substack{\CCC \subset X \\  w(\CCC^{\vphantom{A}} ) < R} } \mathrm{Aire}(\CCC).
\]
Comme pour le comptage des liens de selles, on dispose d'une formule de Siegel-Veech: il existe une constante $C_{\mathrm{aire}}(\nu)$ telle que pour $\nu$-presque toute surface de la strate $\HHH{(k_1, \dots, k_r)}$,
\[
\displaystyle\lim_{T \to +\infty} T^{-2}\,{N_{\mathrm{aire}}(X,T)} = \pi \,C_{\mathrm{aire}}(\nu).
\]

\subsection{Le th\'{e}or\`{e}me d'Eskin-Kontsevich-Zorich}

Eskin, Kontsevich et Zorich \cite{EKZ2} donnent une formule reliant la somme des exposants de Liapounoff \`{a} la constante $C_{\mathrm{aire}}(\nu)$.
Avec les notations pr\'{e}c\'{e}dentes, on a~:

\begin{theorem} \label{cacasselabaraque}
Etant donn\'{e}e une composante connexe $\HHH{(k_1, \dots, k_r)}$ de strate de diff\'{e}rentielles holo\-morphes, soient $\lambda_1, \dots,   \lambda_g$ les exposants de Liapounoff positifs du cocycle de Kontsevich-Zorich associ\'{e}. Alors
\[
\lambda_1 + \dots + \lambda_g = \frac{1}{12} \sum_{i=1}^r \frac{k_i(k_i +2) }{k_i +1} + \frac{\pi^2}{3}C_{\mathrm{aire}} (\nu).
\]
\end{theorem}

\begin{rem}
Une formule analogue existe pour les formes diff\'{e}rentielles quadratiques, nous renvoyons \'{a} l'article original pour plus de d\'{e}tails. Les techniques de preuve sont les m\^{e}mes. 
\end{rem}

Calculer les constantes de Siegel-Veech n'est pas simple. Une fois les formules d'Eskin-Masur \'{e}tablies, il faut comprendre comment calculer explicitement ces constantes. Eskin, Masur et Zorich \cite{EMZ} donnent un algorithme pour les calculer  bas\'{e} sur une \'{e}tude des d\'{e}g\'{e}n\'{e}rescences des strates \`{a} l'infini. La difficult\'{e} est de comprendre comment les sym\'{e}tries des d\'{e}g\'{e}n\'{e}rescences influent sur la valeur des constantes. Le travail d'Eskin-Masur-Zorich  ram\`{e}ne la question au calcul du volume des strates normalis\'{e}es $\HHH_1{(k_1, \dots, k_r)}$, ce qui fait l'objet d'un travail d'Eskin et Okounkov \cite{EO}, lui-m\^{e}me bas\'{e} sur un r\'{e}sultat de Bloch et Okounkov \cite{BO}. La fa\c{c}on de calculer la somme des exposants de Liapounoff du cocycle de Kontsevich-Zorich est donc bien compliqu\'{e}e !
La rationnalit\'{e} de la somme des exposants positifs s'obtient par le m\^{e}me cheminement~: constantes de Siegel-Veech puis volume des strates. Il n'existe pas de preuve plus directe de la rationalit\'{e} connue \`{a} ce jour.

\subsection{Autres mesures invariantes sur les strates}
Les m\^{e}mes questions se posent pour toute mesure de probabilit\'{e} $\textrm{SL}_2(\R)$-invariante sur une sous-orbifolde ferm\'{e}e  $\textrm{SL}_2(\R)$-invariante d'une strate normalis\'{e}e. 
La classification des ferm\'{e}s invariants par  $\textrm{GL}_2(\R)$  est loin d'\^{e}tre achev\'{e}e. N\'{e}anmoins, un pas essentiel vient d'\^{e}tre franchi~:
Eskin et Mirzakhani \cite{EMi} ont annonc\'{e} r\'{e}cemment que tous les ferm\'{e}s $\textrm{GL}_2(\R)$-invariants connexes sont des sous-vari\'{e}t\'{e}s affines dans les coordonn\'{e}es des p\'{e}riodes, 
la mesure associ\'{e}e \'{e}tant proportionnelle \`{a} la mesure de Lebesgue sur le sous-espace correspondant (on parle alors de mesure alg\'{e}brique). 
\par \medskip
Le th\'{e}or\`{e}me d'Eskin-Kontsevich-Zorich est valable pour toute mesure de probabilit\'{e} alg\'{e}brique $\textrm{SL}_2(\R)$-invariante v\'{e}rifiant une condition technique, qui est toujours satisfaite d'apr\`{e}s un r\'esultat d'Avila-Matheus-Yoccoz \cite{AMY}. Ainsi la formule liant  somme des exposants de Liapounoff positifs et constantes de Siegel-Veech est vraie pour toute mesure $\textrm{SL}_2(\R)$-invariante; n\'{e}anmoins la rationalit\'{e} de la somme des exposants positifs est un probl\`{e}me ouvert en toute g\'{e}n\'{e}ralit\'{e}. Ce r\'{e}sultat est connu pour les courbes de Teichm\"{u}ller, c.a.d. pour les orbites ferm\'{e}es sous l'action de $\textrm{SL}_2(\R)$. Dans quelques strates en petit genre, Chen et M\"{o}ller \cite{CM} ont d\'{e}montr\'{e} que 
la valeur de la somme des exposants de Liapounoff positifs est ind\'{e}pendante de la mesure de probabilit\'{e} fix\'{e}e. Ceci devient presque toujours faux d\`{e}s que le genre est assez grand.

\subsubsection{Un exemple simple~: les surfaces \`{a} petits carreaux}

Dans cette section, on s'int\'{e}resse \`{a} l'exemple des surfaces de translation dont l'orbite est ferm\'{e}e. 
Un th\'{e}or\`{e}me de Smillie \cite{SW} affirme que l'orbite de $(X, \omega)$ est ferm\'{e}e si et seulement si le stabilisateur de $(X, \omega)$ est un r\'{e}seau dans $\textrm{SL}_2(\R)$. 
De telles surfaces sont appel\'{e}es surfaces de Veech car celui-ci fut le premier \`{a} les \'{e}tudier syst\`{e}matiquement en 1989\footnote{Il montre par exemple que les polygones r\'{e}guliers avec un nombre pair de c\^{o}t\'{e}s sont des surfaces de Veech lorsqu'on identifie les c\^{o}t\'{e}s parall\`{e}les deux \`{a} deux.}. Le stabilisateur est appel\'{e} le groupe de Veech associ\'{e} \`{a} la surface. 
On appellera courbe de Teichm\"{u}ller $\CCC = \Gamma \backslash\Hyp = \Gamma \backslash \textrm{SL}_2(\R)/\textrm{SO}_2(\R)$ o\`{u} $\Gamma$ est le groupe de Veech de la surface.
 \begin{figure}[h!]
 \begin{center}
 \begin{tikzpicture}
   \tikzstyle{a line}=[very thick, dash pattern=on 8pt off 2pt];
 \draw (0,0)   {}
         -- ++(45:1cm) node [below] {$4$}
         -- ++(45:1cm) 
                 -- ++(90:1cm) node [right] {$1$}
         -- ++(90:1cm) 
         -- ++(135:1cm) node [above] {$2$}
               -- ++(135:1cm)
         -- ++(180:1cm) node [above] {$3$}
                 -- ++(180:1cm) 
         -- ++(225:1cm) node [left] {$4$}
                 -- ++(225:1cm) 
         -- ++(270:1cm) node [left] {$1$}
                 -- ++(270:1cm) node (d1){}
        -- ++(315:1cm) node [below] {$2$}
               -- ++(315:1cm) node (b1){}
                  -- ++(360:1cm)  node [below] {$3$}
         -- (0,0);
         
        \end{tikzpicture}
 \caption{\label{fig:Veechoctogone}
 Une surface de Veech: l'octogone.
 }
 \end{center}
  \end{figure}
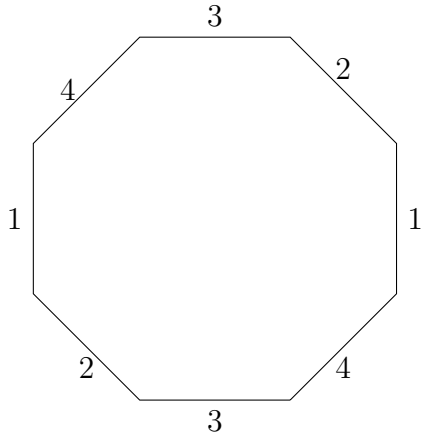 
\par
Les surfaces de Veech ne sont classifi\'{e}es qu'en genre deux par McMullen \cite{Mc1, Mc2, Mc3, Mc4, Mc5, Mc6}, le probl\`{e}me de la classification \'{e}tant totalement ouvert en genre plus grand. 
Il existe cependant une m\'{e}thode simple pour construire de telles surfaces. Le tore $\R^2/\Z^2$ est une surface de Veech car son stabilisateur est $\textrm{SL}_2(\Z)$, qui est un r\'{e}seau de $\textrm{SL}_2(\R)$.
Par d\'{e}finition, une surface \`{a} petits carreaux est un rev\^{e}tement fini du tore  $\R^2/\Z^2$ ramifi\'{e} uniquement au dessus de l'origine \footnote{Les surfaces \`{a} petits carreaux peuvent \'{e}galement s'interpr\'{e}ter comme les points des strates dont les coordonn\'{e}es sont enti\`{e}res.}. De telles surfaces se d\'{e}composent en un nombre fini de carr\'{e}s unit\'{e}s, et elles se d\'{e}composent en cylindres dans toute direction rationnelle.
Il n'est pas tr\`{e}s difficile de montrer qu'une surface \`{a} petits carreaux 
est une surface de Veech, et que son groupe de Veech est commensurable \`{a}  $\textrm{SL}_2(\Z)$; on supposera pour simplifier qu'il s'agit d'un sous-groupe de $\textrm{SL}_2(\Z)$. Le groupe $\textrm{SL}_2(\Z)$ agit sur les surfaces \`{a} petits carreaux par son action lin\'{e}aire, et l'orbite d'une telle surface $X$ est finie; son cardinal est l'indice de son groupe de Veech dans $\textrm{SL}_2(\Z)$. 
Chaque surface $X_i$ de l'orbite de  $\textrm{SL}_2(\Z) \cdot X$ est d\'{e}compos\'{e}e en cylindres horizontaux de p\'{e}rim\`{e}tres $w_{i,j}$ et hauteurs $h_{i,j}$.

 \begin{figure}[h!]
 \begin{center}
 \begin{tikzpicture}
 \tikzstyle{a line}=[very thick, dash pattern=on 8pt off 2pt];
  \draw[very thick] (0,0) -- (4,0);
   \draw[dashed] (0,1) -- (4,1);
   
 \draw[dashed] (1,0) -- (1,2);
  \draw[dashed] (2,0) -- (2,2);
  \draw[dashed] (3,0) -- (3,2);
 
   \draw[very thick] (0,2) -- (4,2);
   \draw[very thick] (0,3) -- (3,3);
 \draw[very thick] (0,4) --(3,4);
 \draw[very thick] (0,0) --(0,4);
 \draw[very thick] (4,0) --(4,2);
 \draw[very thick] (1,2) --(1,3);
 \draw[very thick] (3,3) --(3,4);
 
  \draw[dashed] (1,3) -- (1,4);
  \draw[dashed] (2,3) -- (2,4);
 
  \node at (1.5,0.5) {$\CCC_1$};
 \node at (0.5,2.5) {$\CCC_2$};
 \node at (1.5,3.5) {$\CCC_3$};
  \end{tikzpicture}
 \caption{\label{fig:surface-carreaux}
 Surface \`{a} petits carreaux~: les verticaux c\^{o}t\'{e}s oppos\'{e}s sont identifi\'{e}s. Il y a 3 cylindres horizontaux de param\`{e}tres $(w_1= 4, h_1 = 2)$,  $(w_1= 1, h_2 = 1)$, $(w_3= 3, h_3 = 1)$. 
 }
 \end{center}
 \end{figure}
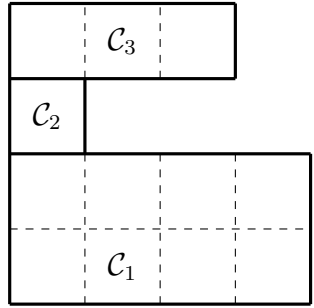
 \par
Les formules pour les constantes de Siegel Veech pour les surfaces \`{a} petits carreaux sont particuli\`{e}rement explicites, ce qui permet d'obtenir le r\'{e}sultat suivant:
\begin{cor}
Si $(X,\omega)$ est une surface \`{a} petits carreaux de la strate $\HHH{(k_1, \dots, k_r)}$, les exposants de Liapounoff du cocycle de Kontsevich-Zorich sont donn\'{e}s par la formule~:
\[
\lambda_1 + \dots + \lambda_g = \frac{1}{12} \sum_{i=1}^r \frac{k_i(k_i +2) }{k_i +1} + \frac{1}{\#\, (\textrm{SL}_2(\Z) \cdot X)} \sum_{X_i \in \textrm{SL}_2(\Z) \cdot X}
 \displaystyle\sum_{\substack {\mathrm{cylindres} \\ \mathrm{horizontaux} \ \CCC_{ij} \\\mathrm{avec} \ X_i = \ \cup  \,\CCC_{ij}}}   \dfrac{h_{ij}}{w_{ij}}\, .
\]
\end{cor}

\subsubsection{Un exemple tr\`{e}s explicite~: les rev\^{e}tements cycliques de la sph\`{e}re}
Forni-Matheus-Zorich \cite{FMZ} et Eskin-Kontsevich-Zorich \cite{EKZ1} ont calcul\'{e} tous les exposants de Liapounoff dans le cas des rev\^{e}tements cycliques de la sph\`{e}re ramifi\'{e}s au-dessus de quatre points. Ces rev\^{e}tements sont obtenus \`{a} partir de courbes affines dans $\mathbb{C}^2$ de la forme
\[
w^N = (z-z_1)^{a_1}(z-z_2)^{a_2}(z-z_3)^{a_3}(z-z_4)^{a_4}
\]
en compl\'{e}tant \`{a} l'infini puis en d\'{e}singularisant. La projection sur la coordonn\'{e}e $z$ d\'{e}finit un rev\^{e}tement cyclique \`{a} valeurs dans la droite projective qui est ramifi\'{e} au dessus des quatre points $z_1$, $z_2$, $z_3$ et $z_4$. Un g\'{e}n\'{e}rateur $T$ du groupe de rev\^{e}tement est donn\'{e} dans ces coordonn\'{e}es par la formule
$T(z,w) = (z, \zeta w)$ o\`{u} $\zeta = \exp(2i\pi/N)$.
Dans le seul cas int\'eressant, le fibr\'{e} de fibre $\textrm{H}^{1,0}$  associ\'{e} \`{a} cette famille de rev\^{e}tements cycliques se scinde en somme directe de fibr\'{e}s en droites holomorphes correspondant aux sous-espaces propres de l'op\'{e}rateur de monodromie $T$. Il s'ensuit que le fibr\'{e} de Hodge r\'{e}el peut \^{e}tre d\'{e}compos\'{e} en sous-fibr\'{e}s invariants par la monodromie (et donc plats) de rang 2 ou 4 dont des sections trivialisantes sur un ouvert de Zariski peuvent \^{e}tre explicitement calcul\'{e}es. Sur chaque sous-fibr\'e un seul exposant de Liapounoff est positif. On peut donc dans ce cas tr\`{e}s particulier calculer {\it tous} les exposants, ce sont des nombres rationnels. En suivant la strat\'{e}gie de la preuve g\'{e}n\'{e}rale du th\'{e}or\`{e}me \ref{cacasselabaraque}, chaque exposant est donn\'{e} comme l'int\'{e}grale de la courbure du fibr\'{e} correspondant. Le point essentiel de cet exemple est que l'espace de param\`{e}tres correspondant aux surfaces est une courbe alg\'{e}brique (orbifolde) compacte. L'int\'{e}grale de courbure s'exprime alors comme le degr\'{e} d'un certain fibr\'{e} en droites orbifolde gr\^{a}ce \`{a} une formule de Peters permettant de contr\^{o}ler les contributions pr\`{e}s des points singuliers. Eskin, Kontsevich et Zorich arrivent \`{a} calculer explicitement ces degr\'{e}s et produisent ainsi les valeurs de tous les exposants dans cette situation. Un cas particuli\`{e}rement fameux est le {\it Eierlegende Wollmilchsau}, un origami \`{a} 8 carreaux qui est un rev\^{e}tement ramifi\'{e} de la sph\`{e}re de genre 3 dont les exposants sont tous nuls sauf le plus grand.
 \begin{figure}[h!]
 \begin{center}
 
 \begin{tikzpicture}
 
 \draw (0,0)   
          -- (1,0) node [below] {$6$}
         -- (2,0) --(2,1) 
                -- (2,2)  
         -- (1,2) node [above] {$4$}
         --(0,2) 
          --(0,1) node [left] {$1$}
           -- (0,0);
           
  \draw (2,0)   
          -- (3,0) node [below] {$7$}
         -- (4,0) --(4,1)
                -- (4,2)  
         -- (3,2) node [above] {$3$}
         --(2,2)
          --(2,1)
           -- (2,0);
          
   \draw (4,0)   
          -- (5,0) node [below] {$8$}
         -- (6,0) --(6,1)
                -- (6,2)  
         -- (5,2) node [above] {$2$}
         --(4,2) 
          --(4,1)
           -- (4,0);      
           
   \draw (6,0)   
          -- (7,0) node [below] {$9$}
         -- (8,0) --(8,1) node [right]{$1$}
                -- (8,2)  
         -- (7,2) 
         --(6,2) 
          --(6,1) 
           -- (6,0);          
              
    \draw (6,2)   
          -- (7,2)  
         -- (8,2)  --(8,3)  
                -- (8,4)  
         -- (7,4) node [above] {$7$}
         --(6,4) 
          --(6,3)  node [left] {$5$}
           -- (6,2);     
           
      \draw (8,2)   
          -- (9,2) node [below] {$2$}
         -- (10,2) --(10,3)
                -- (10,4)  
         -- (9,4) node [above] {$6$}
         --(8,4) 
          --(8,3) 
           -- (8,2);           
           
            \draw (10,2)   
          -- (11,2) node [below] {$3$}
         -- (12,2) --(12,3)
                -- (12,4)  
         -- (11,4)  node [above] {$9$}
         --(10,4) 
          --(10,3)
           -- (10,2);          
           
                \draw (12,2)   
          -- (13,2) node [below] {$4$}
         -- (14,2) --(14,3) node [right] {$5$}
                -- (14,4)  
         -- (13,4) node [above] {$8$}
         --(12,4) 
          --(12,3) 
           -- (12,2);

 \end{tikzpicture}
 
\caption{\label{fig:wollmichsau}
\textit{Eierlegende Wollmichsau} (c\'{e}l\`{e}bre animal l\'{e}gendaire allemand pouvant s'ap\-parenter \`{a} une poule aux \oe{}ufs d'or).}
\end{center}
\end{figure}
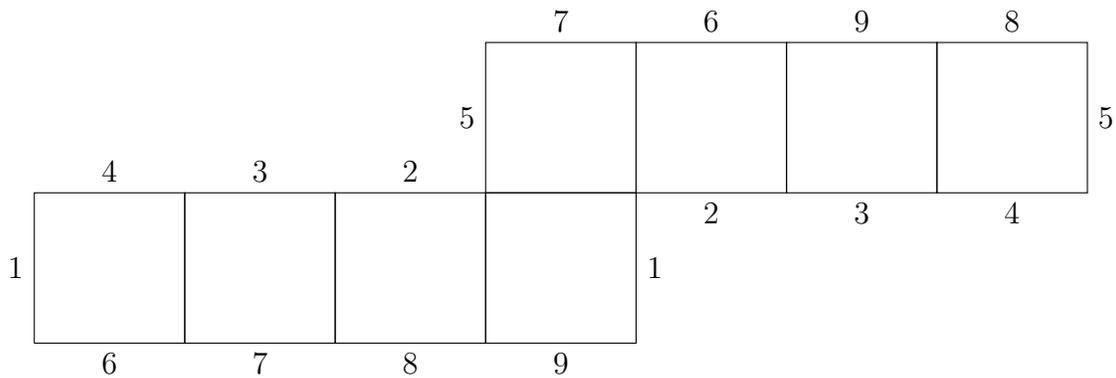

\subsection{Une autre famille d'exemples plus compliqu\'{e}e}
Bouw et M\"{o}ller \cite{BM} ont construit une famille de surfaces de Veech d'une fa\c con nouvelle et non triviale et ont calcul\'{e} explicitement une partie des exposants de Liapounoff du cocycle de Kontsevich-Zorich pour cette famille. La famille contient les exemples de Veech (les polygones r\'{e}guliers). Les m\'{e}thodes de preuves n\'{e}cessitent des arguments alg\'{e}briques suppl\'{e}mentaires que nous ne d\'{e}velopperons pas ici.

\section{Quelques id\'{e}es de preuves}

\subsection{Strat\'{e}gie g\'{e}n\'{e}rale}

Nous donnons ici l'id\'{e}e g\'{e}n\'{e}rale de la preuve, et nous d\'{e}taillerons dans la suite quelques points importants. 
\par\medskip
Dans un premier temps, les auteurs rappellent la formule de Kontsevich \cite{Ko}, d\'{e}montr\'{e}e par Forni \cite{Fo}, qui exprime \`{a} une constante explicite pr\`{e}s la somme des exposants de Liapounoff positifs comme l'int\'{e}grale de la courbure du fibr\'{e} de Hodge dans la direction des disques de Teichm\"{u}ller.
\par\medskip  
Dans le cas des courbes de Teichm\"{u}ller, la famille de courbes peut \^{e}tre compactifi\'{e}e en rajoutant des courbes \'{e}ventuellement nodales aux pointes. La forme de courbure du fibr\'{e} de Hodge est singuli\`{e}re au voisinage des points ajout\'{e}s, cependant un r\'{e}sultat de Peters \cite{Pe} montre que les singularit\'{e}s de cette forme sont mod\'{e}r\'{e}es en un sens d\'{e}fini par Mumford, en particulier elle est int\'{e}grable et son int\'{e}grale est exactement la premi\`{e}re classe de Chern de l'extension de Deligne orbifolde du fibr\'{e} de Hodge \`{a} la courbe compactifi\'{e}e\footnote{La structure orbifolde que l'on met sur la courbe de Teichm\"{u}ller compactifi\'{e}e tient compte de la monodromie de la connexion de Gau\ss-Manin aux pointes.}. Il faut cependant remarquer que cette approche n'implique pas directement un analogue du th\'{e}or\`{e}me \ref{cacasselabaraque}, m\^{e}me si elle prouve la rationalit\'{e} de la somme des exposants et en donne une interpr\'{e}tation cohomologique.
\par \medskip 
On ne sait malheureusement pas construire une compactification qui permette de calculer cohomologiquement l'int\'{e}grale de la courbure du fibr\'{e} de Hodge dans le cas g\'{e}n\'{e}ral. Pour contourner ce probl\`{e}me, on utilise la strat\'{e}gie suivante: comme le fibr\'{e} de Hodge  sur un disque de Teichm\"{u}ller est le pouss\'{e} en avant du fibr\'{e} cotangent relatif de la famille de courbes correspondante, le th\'{e}or\`{e}me de Riemann-Roch analytique permet de calculer la forme de courbure du fibr\'{e} de Hodge. Ce calcul produit un terme de torsion analytique faisant intervenir la fonction z\^{e}ta du Laplacien, qui correspond au terme correctif dans la d\'{e}finition de la m\'{e}trique de Quillen associ\'{e}e.
Les auteurs utilisent un th\'{e}or\`{e}me de Kokotov et Korotkin  \cite{KK}, version ad-hoc du th\'{e}or\`{e}me de Riemann-Roch analytique, qui permet d'exprimer la courbure du fibr\'{e} de Hodge le long des disques de Teichm\"{u}ller en fonction de la fonction z\^{e}ta du laplacien pour la m\'{e}trique plate \footnote{Une autre preuve, bas\'{e}e sur un r\'{e}sultat de Fay, est \'{e}galement d\'{e}velopp\'{e}e dans l'article.}. 
\par \medskip
L'op\'{e}rateur de Laplace pour la m\'{e}trique plate est difficile \`{a} \'{e}tudier, n\'{e}anmoins pr\`{e}s du bord son comportement est proche du laplacien pour la m\'{e}trique hyperbolique\footnote{Cette partie n\'{e}cessite des calculs de g\'{e}om\'{e}trie hyperbolique longs et techniques que nous ne reproduisons pas.}.
Comme le dit Hubbard, une surface de Riemann hyperbolique d\'{e}g\'{e}n\'{e}r\'{e}e est une succession de morceaux de plomberie, reli\'{e}s entre eux par de grands tuyaux. L'image signifie que lorsque la m\'{e}trique hyperbolique d\'{e}g\'{e}n\`{e}re, on dispose \`{a} certains endroits de courbes courtes contenues dans de longs anneaux, et le reste se comporte comme une surface compacte \`{a} bord. 
Les auteurs utilisent un r\'{e}sultat de Lundelius \cite{Lu} qui donne une asymptotique de la fonction z\^{e}ta du Laplacien hyperbolique pour une surface d\'{e}g\'{e}n\'{e}ree en fonction des longueurs hyperboliques courtes. 
Il y a deux types d'anneaux diff\'{e}rents pour la m\'{e}trique plate~: les cylindres plats et les anneaux euclidiens (voir dessins  \ref{fig:anneau} \& \ref{fig:cyl-grand-module}). On montre tout d'abord que les anneaux euclidiens ne contribuent pas dans notre calcul. 
Il faut enfin relier les longueurs des cylindres d\'{e}g\'{e}n\'{e}r\'{e}es aux constantes de Siegel-Veech. 
Ayant des estim\'{e}es fines pr\`{e}s du bord, la preuve consiste \`{a} faire  une sorte d'int\'{e}gration par parties pour ne faire intervenir que le comportement de la fonction z\^{e}ta du laplacien pour la m\'{e}trique plate pr\`{e}s du bord. Cette partie utilise les r\'{e}sultats de Eskin-Masur exprimant les constantes de Siegel-Veech \`{a} l'aide des voisinage du bord des strates. 
\par \medskip
Pour conclure, on peut r\'{e}sumer la strat\'{e}gie de preuve de la mani\`{e}re heuristique suivante: dans le cas des feuilles \textit{alg\'{e}briques} du feuilletage des strates sous l'action de $\textrm{SL}_2(\R)$, c.a.d. des courbes de Teichm\"{u}ller, la situation est naturellement compactifiable et la th\'{e}orie classique des variations de structures de Hodge permet d'interpr\'{e}ter la somme des exposants de Liapounoff de mani\`{e}re cohomologique comme un degr\'{e}. Pour les feuilles transcendantes ou m\^{e}me pour les strates toutes enti\`{e}res, cette m\'{e}thode n'est pour l'instant pas applicable \footnote{M\^{e}me si on peut compactifier les strates, l'interpr\'{e}tation cohomologique de l'int\'{e}grale de la courbure du fibr\'{e} de Hodge n'est pas connue car la base n'est plus une courbe.} et il faut \'{e}tudier les d\'{e}g\'{e}n\'{e}rescences de variations de structures de Hodge au voisinage du bord d'un point de vue \textit{m\'{e}trique}, ce qui est beaucoup plus technique et d\'{e}licat. La rationalit\'{e} provient de calculs purement m\'{e}triques, \`{a} savoir des asymptotiques du nombre de cylindres de grand module plong\'{e}s dans des surfaces d\'{e}g\'{e}n\'{e}r\'{e}es.
\subsection{Formule de Kontsevich}

Cette formule lie la somme des exposants de Liapounoff positifs \`{a} la courbure du fibr\'{e} de Hodge $\mathscr{H}$, qui est le fibr\'{e}  holomorphe de fibre $\mathrm{H}^{1,0}(X)$ au dessus de chaque surface $X$. Elle a \'{e}t\'{e} annonc\'{e}e par Kontsevich puis formellement d\'{e}montr\'{e}e par Forni.
C'est la formule de base de la preuve du th\'{e}or\`{e}me \ref{cacasselabaraque}, nous allons donc l'\'{e}noncer puis donner quelques \'{e}tapes de sa d\'{e}monstration.
\par \bigskip
Soit $\Lambda(X)$ la courbure du fibr\'{e} en droites $\mathscr{H}$ le long des disques de Teichm\"{u}ller d\'{e}finie par 
\[
\Lambda(X) = - \frac{1}{4} \Delta_{\mathrm{Teich}} \log \vert \vert s \vert \vert^2_{\textrm{\tiny{Hodge}}}
\]
o\`{u} $s$ est n'importe quelle section holomorphe locale de $\det \, \mathscr{H}$ qui ne s'annule pas. La quantit\'{e} $\Lambda(X)$ est ind\'{e}pendante de $s$. Si $\omega_1, \dots, \omega_g$ est une base de $1$-formes holomorphes dans un voisinage de $X$, on a
\[
\Lambda(X)= - \frac{1}{4} \Delta_{\mathrm{Teich}} \log \big\vert \det \, \bigl(\iota \, (\omega_i, \omega_j)\bigr) \big\vert.
\]
\begin{theorem} [Kontsevich, Forni] \label{background}

Avec les notations pr\'{e}c\'{e}dentes, on a la relation
\[
\lambda_1 + \dots + \lambda_g = \int_{\HHH_1{(k_1, \dots, k_r)}} \Lambda(X) \, d\nu_1(X).
\]
o\`u $\nu_1$ est la mesure de probabilité obtenue en divisant $\nu$ par son volume.
\end{theorem}
Pour d\'{e}montrer cette identit\'{e}, on utilise le fait bien connu que la somme $\lambda_1 +\dots + \lambda_g$ est l'exposant de Liapounoff maximal pour le cocycle induit sur la $g^{\textrm{e}}$ puissance ext\'{e}rieure $\Lambda^g \,{\mathrm{H}}^1(X,\R)$. Notons $\LLL$ la grasmannienne relative des sous-espaces lagrangiens pour le fibr\'{e} de Hodge et fixons 
une mesure $\sigma$ sur $\LLL$ donn\'{e}e par une forme volume. La premi\`{e}re \'{e}tape est de remarquer que pour $\nu_1$ presque tout $(X, \omega)$ et pour $\sigma$ presque tout $L$ dans $\LLL_X$, si $L_t$ d\'{e}signe le transport de $L$ par la connexion de Gau\ss-Manin le long du flot g\'{e}od\'{e}sique\footnote{Comme la forme d'intersection est plate, les $L_t$ restent lagrangiens.}, 
\begin{equation}\label{eq:somme-exposants}
\lambda_1 + \dots + \lambda_g = \lim_{T \to +\infty} T^{-1}\, \log \vert \vert L_t \vert\vert_{g_tX}.
\end{equation}
D'apr\`{e}s le th\'{e}or\`{e}me d'Oseledets, cette \'{e}quation est satisfaite d\`{e}s que $L$ n'est pas contenu dans le sous-espace vectoriel $H(\omega)$ associ\'{e} aux exposants de Liapounoff autres que l'exposant dominant, \`{a} savoir $E_2(\omega) \oplus \dots \oplus E_{2g}(\omega)$ avec les notations de \S\ref{exposants}. 
Or il est bien connu que pour un cocycle symplectique, le sous-espace d'Oseledets associ\'{e} aux exposants positifs $\lambda_1, \dots, \lambda_g$ est lagrangien. Ainsi le sous-espace $E_1(\omega) \oplus \dots \oplus E_g(\omega)$ est un \'{e}l\'{e}ment de $\LLL_X$ non contenu dans $H(\omega)$. Par cons\'{e}quent, pour presque tout $(X, \omega)$, l'ensemble des sous-espace langrangiens contenu dans $H(\omega)$ est de mesure nulle vu que $\LLL_X$  est une vari\'{e}t\'{e} alg\'{e}brique irr\'{e}ductible. On obtient ainsi \eqref{eq:somme-exposants}.
\par \medskip
Notons $P\mathcal{H}(k_1, \dots, k_r)$, la strate projectivis\'ee obtenue en identifiant deux 1--formes multiples l'une de l'autre.
L'\'{e}tape suivante consiste \`{a} int\'{e}grer la fonction constante $\displaystyle\lim_{T \to \infty}T^{-1}\, \log \vert \vert L_t \vert\vert_{g_tX}$ sur la strate projectivis\'ee $P\mathcal{H}(k_1, \dots, k_r)$, puis \`{a} effectuer des manipulations formelles (formule de Green, interversion de limites) dues \`{a} Forni dont certaines n\'{e}cessitent des justifications non triviales. Expliquons o\`{u} l'utilisation de sous-espaces \textit{lagrangiens} intervient. 
\par \medskip
Comme les strates sont feuillet\'{e}es par les $\textrm{GL}_2(\R)$-orbites des surfaces de translation, les strates projectivis\'ees sont feuillet\'{e}es par des disques de Teichm\"{u}ller, c'est-\`{a}-dire par des quotients du plan hyperbolique $\textrm{SL}_2(\R) / \textrm{SO}_2(\R)$. On munit ces disques de la m\'{e}trique hyperbolique de courbure $-4$, ce qui est coh\'{e}rent avec notre d\'{e}finition du flot g\'{e}od\'{e}sique introduite pr\'{e}c\'{e}demment, et on note $\Delta_{\mathrm{Teich}}$ le laplacien hyperbolique associ\'{e} sur la strate projectivis\'ee. Fixons une base de sections holomorphes locales $\omega_1, \dots, \omega_g$ au voisinage de $X$ ainsi que des bases $v_1(t), \dots, v_g(t)$ de $L_t$ localement constantes en $t$. Alors la quantit\'{e} $\Delta_{\mathrm{Teich}} \log \vert \vert v_1(t) \wedge \dots \wedge v_g(t) \vert \vert _{\textrm{\tiny{Hodge}}}$ ne d\'{e}pend pas du choix de $v_1(t), \dots, v_g(t)$, on la note $\Delta_{\mathrm{Teich}} \log \vert \vert L \vert \vert _{\textrm{\tiny{Hodge}}}$. Un point beaucoup plus remarquable est que cette quantit\'{e} ne d\'{e}pend pas non plus du sous-espace lagrangien initial $L$ en vertu de l'identit\'{e}
\begin{equation} \label{eq:deltaL}
\Delta_{\mathrm{Teich}} \log \vert \vert L \vert \vert _{\textrm{\tiny{Hodge}}}= - \frac{1}{2} \Delta_{\mathrm{Teich}} \log \big\vert \det\, \bigl(\iota \, ( \omega_i, \omega_j )\bigr)_{i,j} \big\vert \, .
\end{equation} 
La preuve, qui ne fait intervenir que de l'alg\`{e}bre multilin\'{e}aire \'{e}l\'{e}mentaire, est la suivante: si $v_1, \dots, v_g$ est une base de $L$, le carr\'{e} de la norme de Hodge de $v_1 \wedge \dots \wedge v_g$ est le d\'{e}terminant de la matrice de  Gram r\'{e}elle associ\'{e} aux vecteurs $v_1, \dots, v_g$, soit en utilisant les notations de \S  \ref{subsection:Hodge}
\[
\vert \vert v_1 \wedge \dots \wedge v_g\vert \vert ^2_{\textrm{\tiny{Hodge}}} = \det\, \bigl( \mathrm{Re} \, \iota \,\bigl(\phi^{-1} (v_i), \phi^{-1} (v_j)\bigr) \bigr)_{1 \leq i,j \leq g}.
\] 
Comme $L$ est Lagrangien, les coefficients $\iota  \,\bigl(\phi^{-1} (v_i), \phi^{-1} (v_j)\bigr)$ sont r\'{e}els, donc $\vert \vert v_1 \wedge \dots \wedge v_g\vert \vert ^2_{\textrm{\tiny{Hodge}}} $ est \'{e}gal au d\'{e}terminant de la matrice de Gram \textit{complexe} $\bigl(\iota \,\bigl(\phi^{-1} (v_i), \phi^{-1} (v_j)\bigr) \bigr)_{1 \leq i,j \leq g}$. Si on munit  $\Lambda_{\C}^g \,\mathrm{H}^{1,0}(X)$ de la norme induite par $\iota$, on a donc 
\[
\vert \vert v_1 \wedge \dots \wedge v_g\vert \vert ^2_{\textrm{\tiny{Hodge}}}=\vert \vert \phi^{-1}(v_1) \wedge \dots \wedge \phi^{-1} (v_g) \vert \vert _{\Lambda_{\C}^g \mathrm{H}^{1,0}(X)}.
\]
Notons $\Omega $ la forme volume sur $\Lambda^{2g}_{\C} \,{\mathrm{H}}^1(X,\C)$ donn\'{e}e par la $g^{\textrm{e}}$ puissance exterieure de la forme symplectique et posons 
\[
\vert \vert v_1 \wedge \dots \wedge v_g \vert \vert^2_{\textrm{\tiny{EKZ}}} = \displaystyle \frac{ \vert \Omega(v_1 \wedge \dots \wedge v_g \wedge \omega_1 \wedge \dots \wedge \omega_g)\vert ^2}{  \vert \Omega(\omega_1 \wedge \dots \wedge \omega_g \wedge \overline{\omega}_1 \wedge \dots \wedge \overline{\omega}_g) \vert} \, .
\]
La quantit\'{e} $\vert \vert \, . \, \vert \vert_{\textrm{\tiny{EKZ}}}$ d\'{e}finit apr\`{e}s transport par $\phi$ une norme sur 
$\Lambda_{\C}^g \mathrm{H}^{1,0}(X)$\footnote{Il est important de remarquer que $\vert\vert.\vert\vert_{\textrm{\tiny{EKZ}}}$ n'est pas d\'{e}finie positive sur $\Lambda_{\R}^g \mathrm{H}^{1}(X, \R)$, ce n'est donc pas une norme mais seulement une pseudo-norme. Ceci montre que la formule $(\ref{eq:formuleKontsevich})$ ci-apr\`{e}s n'est pas valable en g\'{e}n\'{e}ral pour des vecteurs $v_1, \dots, v_g$ arbitraires.}. Comme $\Lambda_{\C}^g \mathrm{H}^{1,0}(X)$ est de dimension complexe $1$, il existe une constante $k$ ind\'{e}pendant des $v_i$ telle que 
\begin{equation}  \label{eq:formuleKontsevich}
\vert \vert v_1 \wedge \dots \wedge v_g \vert \vert^2_{\textrm{\tiny{Hodge}}}=k\, \vert \vert v_1 \wedge \dots \wedge v_g \vert \vert^2_{\textrm{\tiny{EKZ}}}
\end{equation}
et on v\'{e}rifie facilement que $k=1$.
La formule $(\ref{eq:deltaL})$ d\'{e}coule alors directement de $(\ref{eq:formuleKontsevich})$.
Remarquons pour cl\^{o}turer cette discussion que le terme de courbure $\Lambda(X)$ appara\^{\i}t naturellement dans la formule $(\ref{eq:deltaL})$. 
\par \bigskip
Lorsque $\CCC$ est une courbe de Teichm\"{u}ller, le Th\'{e}or\`{e}me devient
\begin{equation} \label{eq:Kontsevich-courbes}
\lambda_1 + \dots + \lambda_g = \dfrac{{i} \displaystyle\int_\CCC \Theta}{\pi (2 \mathfrak{g}_{\CCC} -2 + s_{\CCC})}
\end{equation}
o\`{u} $\Theta$ est la $2$-forme de courbure du fibr\'{e} $\det {\mathscr{H}}$, $\mathfrak{g}_{\CCC}$ est le genre de la coube de Teichm\"{u}ller \footnote{Le genre de la courbe de Teichm\"{u}ller n'a a priori rien \`{a} voir avec celui de la surface de Veech $X$ qui l'engendre.} et $s_{\CCC}$ est le nombre de pointes de $\CCC$. 
Dans cette formule le d\'{e}nominateur est la caract\'{e}ristique d'Euler Poincar\'{e}, c'est-\`{a}-dire l'aire de la surface hyperbolique. Ceci provient du fait que la formule de Kontsevich est obtenue avec une mesure de probabilit\'{e}. Dans ce cas particulier, on peut exprimer l'int\'{e}grale de la courbure comme un degr\'{e} de fibr\'{e}, gr\^{a}ce \`{a} la formule de Peters
\[
\frac{i}{2\pi} \int_{\CCC} \Theta = \deg \overline{\mathscr{H}}
\] 
o\`{u} $\overline{\mathscr{H}}$ est l'extension de Deligne orbifolde du fibr\'{e} de Hodge aux pointes pour une structure orbifolde sur $\overline{\CCC}$ permettant de rendre unipotente la monodromie de la connexion de Gau\ss-Manin aux pointes. C'est gr\^{a}ce \`{a} cette formule que l'on peut calculer tous les exposants de Liapounoff pour les rev\^{e}tements cycliques de la sph\`{e}re. En effet d\`{e}s que l'on dispose d'une sous-structure de Hodge invariante par notre cocycle,  la formule de Kontsevich et la formule de Peters s'appliquent \textit{mutatis mutandis}. C'est \'{e}galement l'architecture de preuve utilis\'{e}e par Bouw et M\"{o}ller pour calculer des exposants de Liapounoff de courbes de Teichm\"{u}ller. Dans ces derniers exemples, la strat\'{e}gie pour casser le fibr\'{e} de Hodge en morceaux invariants est d'utiliser le fait que le plan engendr\'{e} par $\mathrm{Re}(\omega)$ et  $\mathrm{Im}(\omega)$ est d\'{e}fini sur un corps de nombres, et de faire intervenir l'action du groupe de Galois.

\subsection{D\'{e}terminant du Laplacien et formule de Riemann-Roch analytique}
Soit $X$ une surface de Riemann compacte, $g$ une m\'{e}trique riemannienne lisse sur $X$, $\Delta_g$ le laplacien associ\'{e} \footnote{Nous prenons la convention de signe des analystes pour le laplacien: toutes ses valeurs propres sont n\'{e}gatives.} et $(\mu_i)_{i \geq 0}$ le spectre de $-\Delta_g$. La fonction z\^{e}ta du laplacien $\Delta_g$ est d\'{e}finie pour $Re(s)>1$ par la formule $\zeta(s) = \sum_{i=0}^\infty \mu_i^{-s}$. Les estimations classiques du noyau de la chaleur entra\^{\i}nent l'existence d'un prolongement analytique de la fonction z\^{e}ta en une fonction m\'{e}romorphe sur $\C$ avec $1$ pour seul p\^{o}le. Le d\'{e}terminant du laplacien est d\'{e}fini par, $\det \Delta_g = \exp(- \zeta'(0))$, et correspond formellement au produit infini des $\mu_i$ (qui est un produit infini divergent). Pour de plus amples d\'{e}tails, on renvoie le lecteur aux r\'{e}f\'{e}rences \cite[Chap. 2]{BGV} et \cite[Chap. V]{So}.
\par \bigskip

Comme expliqu\'{e} pr\'{e}c\'{e}demment, il va \^{e}tre n\'{e}cessaire de comparer la m\'{e}trique hyperbolique et la m\'{e}trique plate. Pour ce faire, rappelons la formule de Polyakov. Soient $g_1$ et $g_2$ deux m\'{e}triques riemannienes lisses dans la m\^{e}me classe conforme, et \'{e}crivons $g_2 = \exp(2\phi) g_1$ o\`{u} $\phi$ est une fonction lisse \`{a} valeurs r\'{e}elles. Alors on a:

\begin{align} \label{eq:Polyakov}
\log \det \Delta_{g_2} - \log \det \Delta_{g_1} = \frac{1}{12 \pi}\Bigl( \int_X \phi \Delta_{g_1} \phi \, dg_1 - 2 \int_X &\phi K_{g_1} dg_1\Bigr)\\
&+\log \mathrm{Aire}_{g_2} (X) - \log \mathrm{Aire}_{g_1}(X). \notag
\end{align}

La m\'{e}trique plate n'est pas lisse aux singularit\'{e}s de la surface de translation, on peut n\'{e}anmoins d\'{e}finir un d\'{e}terminant du laplacien relatif. Fixons une surface de translation $(X_0, \omega_0)$ dans une strate $\HHH(k_1, \dots, k_r)$, et faisons varier $(X, \omega)$ dans la m\^{e}me strate. On pose
\[
\det \Delta_{\mathrm{plat}}(X, X_0) = \lim_{\varepsilon \to 0} \, \frac{ \det \Delta_{\mathrm{plat}, \varepsilon}(X)}{ \det \Delta_{\mathrm{plat}, \varepsilon}(X_0)}
\]
o\`{u} $ \Delta_{\mathrm{plat}, \varepsilon}(X)$ est une approximation de la m\'{e}trique plate au niveau des z\'{e}ros (on ne modifie la m\'{e}trique que dans un $\varepsilon$ voisinage des z\'{e}ros par une fonction radiale).
\par\medskip
Sur la surface $(X, \omega)$, il y a une unique m\'{e}trique hyperbolique dans la classe conforme de la m\'{e}trique plate poss\'{e}dant des pointes aux z\'{e}ros de la 1--forme $\omega$. La surface de Riemann associ\'{e}e n'est pas compacte, il faut donc aussi r\'{e}gulariser la m\'{e}trique pour d\'{e}finir le d\'{e}terminant du laplacien, ce qui conduit \`{a} nouveau \`{a} un d\'{e}terminant du laplacien relatif $\det \Delta_{\mathrm{hyp}}(X, X_0)$.

\begin{theorem} \label{thm:RRanalytique}
Soit $X_0$ une surface fix\'{e}e dans la strate $\HHH(k_1 \dots k_r)$. Alors, pour tout $(X, \omega)$ dans cette strate,
\[
-4 \Lambda(X) = \Delta_{\mathrm{Teich}} \log \det \Delta_{\mathrm{plat}}(X, X_0) - \frac{1}{3} \sum_{j=1}^r \frac{k_j(k_j + 2)}{k_j +1} \,.
\]
\end{theorem}

Le terme combinatoire est celui du th\'{e}or\`{e}me principal, il restera  donc \`{a} relier $\Delta_{\mathrm{Teich}} \log \det \Delta_{\mathrm{plat}}(X, X_0)$ aux constantes de Siegel-Veech.
\par \medskip
Expliquons comment v\'{e}rifier la formule dans le cas du tore ce qui est d\'{e}j\`{a} fort instructif. Pour ce faire, 
on consid\`{e}re la coordonn\'{e}e $\tau$ dans le domaine fondamental de la surface modulaire 
\[
\mathrm{Im}(\tau) >0, \quad\vert \tau \vert \geq 1, \quad-\frac{1}{2} \leq \mathrm{Re} (\tau) < \frac{1}{2} \, .
\]
Pour tout tore d'aire $1$, la formule de Ray-Singer \cite{RS} s'\'{e}crit
\[
\det\Delta_{\mathrm{plat}} = 4\, \mathrm{Im} (\tau) \,\vert \eta(\tau) \vert^4
\]
o\`{u} $\eta$ est une forme modulaire explicite.
Ainsi 
\[
\Delta_{\mathrm{Teich}} \log \det\Delta_{\mathrm{plat}} = \Delta_{\mathrm{Teich}} \log\,  \mathrm{Im}  (\tau)
\] 
car $\eta$ est holomorphe. 
On peut bien s\^{u}r calculer dans ce cas $\big\vert \det \, \bigl(\iota \, (\omega_i, \omega_j)\bigr) \big\vert$: il y a une seule forme holomorphe $\omega$ de p\'{e}riodes 1 et $\tau$, de sorte que
$\big\vert \det \, \bigl(\iota \, (\omega_i, \omega_j)\bigr) \big\vert= \vert\vert \omega \vert \vert^2 = \mathrm{Im} (\tau)$ ce qui \'{e}tablit la formule souhait\'{e}e dans ce cas tr\`{e}s simple. 
\par \medskip
Une mani\`{e}re de prouver le th\'{e}or\`{e}me  \ref{thm:RRanalytique} est d'utiliser un th\'{e}or\`{e}me de Kokotov et Korotkin \cite{KK} dont l'\'{e}nonc\'{e} est le suivant:
\begin{theorem} [Kokotov, Korotkin]
Pour toute surface de translation $(X, \omega)$ et tout point fix\'{e} $X_0$ dans la strate $\HHH(k_1, \dots, k_r)$, on a 
\[
\det \Delta_{\mathrm{plat}}(X, X_0) = \mathrm{k} \cdot \, \mathrm{Aire}(X, \omega)\, \det\, ( \mathrm{Im} \, B) \cdot \vert \tau(X, \omega)\vert^2
\]
o\`{u} $B$ est la matrice des p\'{e}riodes secondaires de $X$ et $\tau$ est une section holomorphe
d'un fibr\'{e} en droite au-dessus de la strate $\HHH(k_1, \dots, k_r)$.
De plus, 
$\tau(X,\omega)$ est homog\`{e}ne en $\omega$ de degr\'{e} 
\[
p = \frac{1}{12} \sum_{i=1}^r \frac{k_i (k_i +2)}{k_i+1} \, .
\]
\end{theorem}

On voit clairement que ce th\'{e}or\`{e}me est une g\'{e}n\'{e}ralisation de la formule de Ray-Singer. 
\par \bigskip
Comme nous l'avons expliqu\'{e} auparavant, le th\'{e}or\`{e}me \ref{thm:RRanalytique} est un r\'{e}sultat de type Riemann-Roch analy\-tique. Pour motiver ce th\'{e}or\`{e}me, nous allons \'{e}tablir un r\'{e}sultat plus faible (car cohomologique) en appliquant le th\'{e}or\`{e}me de Grothendieck-Riemann-Roch usuel pour des familles de courbes. Ceci permettra d'expliquer de mani\`{e}re cohomologique l'apparition du myst\'{e}rieux terme combinatoire $\dfrac{1}{12} \displaystyle \sum_{i=1}^r \frac{k_i (k_i +2^{\vphantom{A^{\vphantom{A)}}}})}{k_i+1}\,.$
\par \medskip
Soit $\pi \colon \mathcal{C} \rightarrow B$ une famille holomorphe de surfaces de translation dans la strate $\mathcal{H}(k_1, \dots, k_r)$ param\'{e}tr\'{e}e par une base $B$. On dispose d'une section holomorphe globale du fibr\'{e} projectif $\mathbb{P}(\pi_*\Omega^1_{\mathcal{C}/B})$ sur $B$. Cette section peut \^{e}tre interpr\'{e}t\'{e}e comme une section globale $\omega$ partout non nulle de $\pi_* \Omega^1_{\mathcal{C}/B} \otimes \mathcal{L}^{\vee}$, o\`{u} $\mathcal{L}$ est un fibr\'{e} en droites holomorphe sur $B$. Quitte \`{a} prendre un rev\^{e}tement fini sur la base, il existe $r$ sections $s_1, \dots, s_r$ de $\mathcal{C}$ telles que si $D_1, \dots, D_r$ sont les images de ces sections, le diviseur d'annulation de $\omega$ (en tant que section globale sur $\mathcal{C}$ de $\Omega^1_{\mathcal{C}/B} \otimes \pi^*\mathcal{L}^{\vee}$) est $k_1 D_1 + \dots + k_r D_r$. Cela signifie concr\`{e}tement que chaque forme $\omega_b$ d\'{e}finissant la structure de translation de $\mathcal{C}_b$ s'annule aux points marqu\'{e}s $s_1(b), \dots, s_r(b)$ avec multiplicit\'{e}s respectives $k_1, \ldots, k_r$. On a donc l'isomorphisme
\[
\Omega^1_{\mathcal{C}/B} \simeq \pi^* \mathcal{L} \otimes \mathcal{O}\Bigl(\, \sum_{i=1}^r \, k_i D_i\, \Bigr).
\]  
\begin{prop} Si $\mathscr{H}$ est le fibr\'{e} de Hodge $\pi_* \Omega^1_{\mathcal{C}/B}$ sur $B$,
\[
\mathrm{c}_1(\mathscr{H})=\frac{1}{12} \sum_{i=1}^r \frac{k_i (k_i +2)}{k_i+1} \,\mathrm{c}_1(\mathcal{L}).
\]
\end{prop}
Pour prouver ce r\'{e}sultat, on applique le th\'{e}or\`{e}me de Grothendieck-Riemann-Roch pour le morphisme $\pi$ au faisceau $\mathcal{O}_{\mathcal{C}}$. On a $\mathrm{R}^0 \pi_* \mathcal{O}_{\mathcal{C}}=\mathcal{O}_B$ et en utilisant la dualit\'{e} de Grothendieck relative, 
\[
\mathrm{R}^1 \pi_* \mathcal{O}_{\mathcal{C}} \simeq (\mathrm{R}^0 \pi_* \Omega^1_{\mathcal{C}/B})^{\vee}  \simeq \mathscr{H}^{\vee}
\]
donc la classe de $\mathrm{R} \pi_* \mathcal{O}_{\mathcal{C}}$ dans $\mathrm{K}(B)$ est $[\mathcal{O}_B]-[\mathscr{H}^{\vee}]$. On en d\'{e}duit l'\'{e}galit\'{e} \[
\mathrm{c}_1(\mathscr{H})=\pi_*\bigl(\mathrm{td}_2(\mathrm{T}_{\mathcal{C}/B})\bigr)=\frac{1}{12} \pi_* \bigl(\mathrm{c}_1(\mathrm{T}_{\mathcal{C}/B})^2 \bigr) \, .
\]
Comme $\mathrm{c}_1(\mathrm{T}_{\mathcal{C}/B})=-\pi^* \mathrm{c}_1(\mathcal{L})-\sum_{i=1}^r k_i [D_i]$, 
\[
\pi_* \bigl(\mathrm{c}_1(\mathrm{T}_{\mathcal{C}/B})^2 \bigr)=2\sum_{i=1}^r k_i\,\mathrm{c}_1(\mathcal{L}) \, \pi_* [D_i] + \sum_{i=1}^r k_i^2\, \pi_*[D_i]^2.
\]
Remarquons que $\pi_*\, [D_i]=\pi_* {s_i}_*(1)=1$. De plus, la suite exacte des diff\'{e}rentielles relatives restreinte \`{a} $D_i$ entra\^{\i}ne l'\'{e}galit\'{e} $(\Omega^1_{\mathcal{C}/B})_{|D_i}=N^*_{D_i/X}$, donc en identifiant $D_i$ et $B$, $N_{D_i/X}^{\, \otimes k_i+1}=\mathcal{L}^{\vee}$.
On en d\'{e}duit 
\[
\pi_* \, [D_i]^2=\pi_* ({s_i}_*(1). [D_i])=\pi_* {s_i}_* {s_i}^* [D_i]=\mathrm{c}_1(N_{D_i/X})=-\frac{\mathrm{c}_1(\mathcal{L})}{k_i+1}
\]
d'o\`{u} la formule souhait\'{e}e
\[
\pi_* \bigl(\mathrm{c}_1(\mathrm{T}_{\mathcal{C}/B})^2 \bigr)=\sum_{i=1}^r \Bigl(2k_i-\frac{k_i^2}{k_i+1}\Bigr)\,  \mathrm{c_1}(\mathcal{L})=\sum_{i=1}^r \frac{k_i (k_i +2)}{k_i+1} \,\mathrm{c}_1(\mathcal{L}).
\]

\subsection{Compactification et th\'{e}or\`{e}mes de Rafi}

On veut comprendre ce qui se passe lorsque les surfaces de Riemann d\'{e}g\'{e}n\`{e}rent, il faut donc d\'{e}finir une compactification adapt\'{e}e \`{a} la situation. On peut faire cela en termes alg\'{e}briques, ce qui correspond \`{a} la compactification classique de Deligne-Mumford. Topologiquement, on rajoute pour compactifier des surfaces de Riemann dont certaines courbes ont \'{e}t\'{e} pinc\'{e}es. Il est cependant difficile de comprendre alg\'{e}briquement la g\'{e}om\'{e}trie plate \`{a} l'infini. Nous allons maintenant d\'{e}velopper un point de vue topologique d\^{u}e \`{a} Rafi (\cite{Ra1}, \cite{Ra2}) permettant de bien exprimer la g\'{e}om\'{e}trie plate dans la compactification.
\par
\begin{figure}[ht]
 \begin{center}
 \begin{tikzpicture}[scale=2]

 \draw[ arrows=<-]  (0,-0.3)--(1.5,-0.3) node[above] {$1$};
 \draw[ arrows=->]  (1.5,-0.3)--(3,-0.3);
 \draw[ arrows=<-]  (-0.3,0)--(-0.3,1.5) node[right] {$1$};
 \draw[ arrows=->]  (-0.3,1.5)--(-0.3,3);
 
 \draw[very thick]  (0,0)rectangle(3,3);
 \draw[very thick]  (4,0)rectangle(4.8,0.8);
 \draw[very thick]  (1.5,1.2)--(1.5,1.5) node[left]{$a$ } node[right]{$b$ }--(1.5,1.8);
 \draw[very thick]  (4.4,0.1)--(4.4,0.4)node[left]{$b$ } node[right]{$a$ }--(4.4,0.7);
 
 \draw[arrows=<-]  (4,-0.3)--(4.4,-0.3) node[above] {$\varepsilon$};
 \draw[arrows=->]  (4.4,-0.3)--(4.8,-0.3);
 \draw[arrows=<-]  (3.7,0)--(3.7,0.4) node[right] {$\varepsilon$};
 \draw[arrows=->]  (3.7,0.4)--(3.7,0.8);

 \end{tikzpicture}
 \end{center}
  \caption{Deux tores recoll\'{e}s le long d'un lien de selles, celui de droite est d'aire beaucoup plus petite que celui de gauche.  }
 \label{fig:2tores}
  \end{figure}
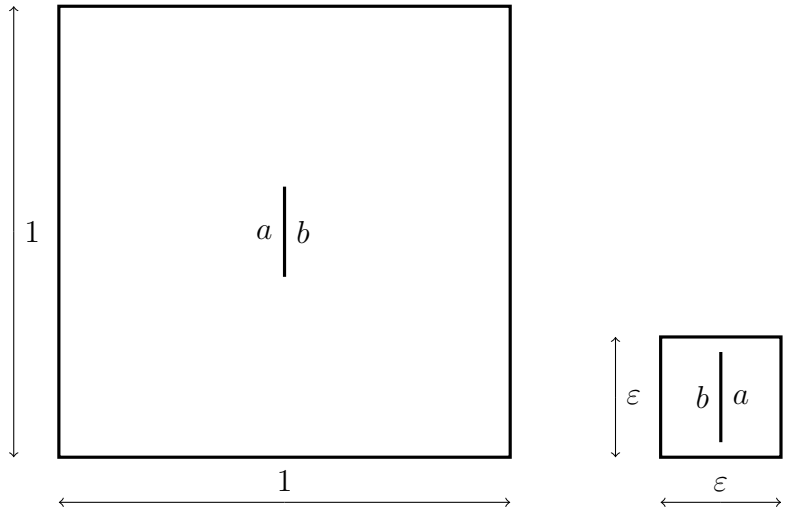

 \medskip
Tout d'abord, on rappelle que la surface de translation $(X, \omega)$ est munie de la m\'{e}trique plate associ\'{e}e mais aussi de la m\'{e}trique hyperbolique se trouvant dans la m\^{e}me classe conforme et poss\'{e}dant des pointes aux z\'{e}ros de la $1$--forme holomorphe $\omega$. Hors d'un compact de l'espace des modules des surfaces de Riemann de genre $g$, on sait exactement localiser les g\'{e}od\'{e}siques courtes. Deux g\'{e}od\'{e}siques {\it hyperboliques} courtes ne peuvent s'intersecter donc leur nombre est born\'{e} par $3g-3+r$ o\`{u} $g$ est le genre et $r$ le nombre de z\'{e}ros.
Pour $\delta$ assez petit, notons $\Gamma(\delta)$ l'ensemble des g\'{e}od\'{e}siques de longueur au plus $\delta$. 
Sur une surface de Riemann de volume fini, une courbe $\delta$--courte est l'\^{a}me d'un anneau de grand module. Il faut par contre noter que de telles courbes ne sont en g\'{e}n\'{e}ral pas g\'{e}od\'{e}siques pour la m\'{e}trique plate.
Au niveau plat, il y a deux possibilit\'{e}s, d\'{e}crites sur les figures \ref{fig:anneau} \& \ref{fig:cyl-grand-module}: cet anneau est ou bien un cylindre plat, ou bien un anneau euclidien.
  \begin{figure}[h!]
  \begin{center}
  \begin{tikzpicture}[scale=2] 
  \fill[fill=black!20] (1.5,1.5) circle (1.5);
  \fill[white] (1.5,1.5) circle (0.3); 
 \draw[very thick]  (0,0)rectangle(3,3);
 \draw[very thick]  (5,0)rectangle(8,3);
 \draw[very thick]  (1.5,1.2)--(1.5,1.5) node[left]{$a$ } node[right]{$b$ }--(1.5,1.8);
 \draw[very thick]  (6.5,1.2)--(6.5,1.5)node[left]{$b$ } node[right]{$a$ }--(6.5,1.8);
 \draw[very thick]  (1.5,1.5) circle (0.6); 
 \end{tikzpicture}
 \end{center}
 \caption{Deux tores sont recoll\'{e}s le long d'un lien de selles. Les c\^{o}t\'{e}s des carr\'{e}s sont identifi\'{e}s, les lien de selles sont identifi\'{e}s suivant $a$ et $b$. Le lien de selles est beaucoup plus court que le c\^{o}t\'{e} du carr\'{e}, la courbe repr\'{e}sent\'{e}e par un cercle est courte. On voit en gris l'anneau de grand module qui l'entoure.  }
 \label{fig:anneau}
  \end{figure}
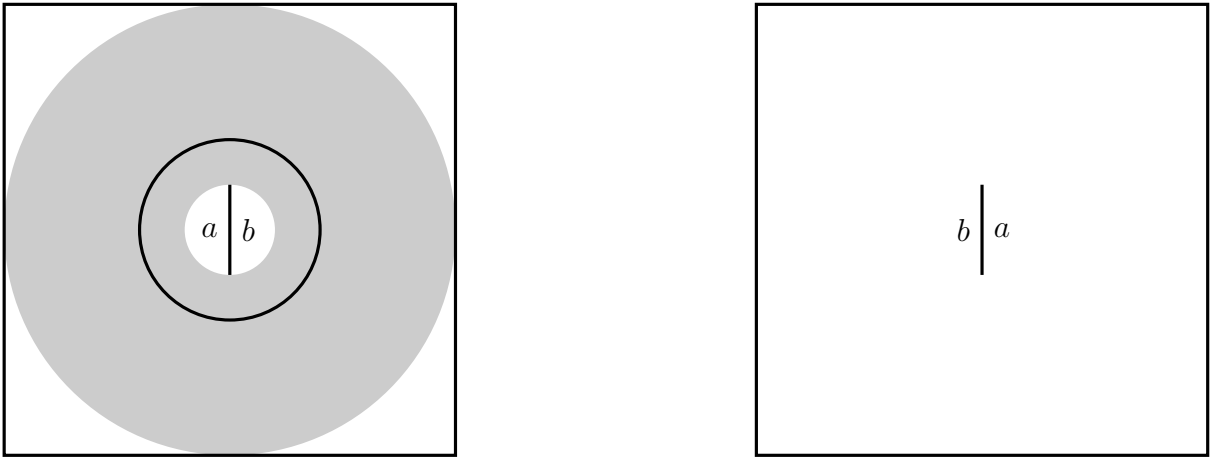 
  \begin{figure}[h!]
  \begin{center}
  \begin{tikzpicture}[scale=1]
  \draw[very thick]  (0,0)rectangle(1,5);
  \draw[dashed]  (0,2.5)--(1,2.5);
  \end{tikzpicture}
  \end{center}
  \caption{Cylindre plat de grand module, les c\^{o}t\'{e}s verticaux sont identifi\'{e}s. La courbe en pointill\'{e} est courte.}
  \label{fig:cyl-grand-module}
  \end{figure}
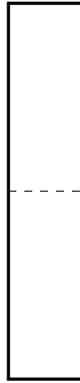
Pour  $\gamma$ courbe simple ferm\'{e}e non p\'eriph\'erique, on consid\`{e}re son repr\'{e}sentant g\'{e}od\'{e}sique pour la m\'{e}trique {\it plate}, il est unique sauf si $\gamma$ est l'\^{a}me d'un cylindre. On d\'{e}coupe alors la surface $X$ le long des repr\'{e}sentants g\'{e}od\'{e}siques de $\Gamma(\delta)$, et on \^{o}te compl\`{e}tement les cylindres plats dont l'\^{a}me est dans $\Gamma(\delta)$.
Une composante connexe du complementaire de ces cylindres est appel\'{e}e $\delta$--\'{e}paisse. Etant donn\'{e} une telle composante $Y$, on obtient en d\'{e}coupant le long des repr\'{e}sentants plats de $\Gamma(\delta)$ une sous-surface \`{a} bord $\hat{Y}$, appel\'{e}e par Rafi repr\'{e}sentant plat de $Y$, qui a le m\^{e}me type d'homotopie que $Y$. Il y a ici  une subtilit\'{e} que nous passerons sous silence: la sous-surface $\hat{Y}$ peut d\'{e}g\'{e}n\'{e}rer sous forme de graphe (il faut alors consid\'{e}rer qu'elle vient avec un voisinage infinit\'{e}simal pour que le discours reste correct). 
Rafi \cite{Ra1}, \cite{Ra2} d\'{e}finit la taille $\lambda(\hat{Y})$ d'une sous-surface \`{a} bord Y comme \'{e}tant la longueur plate minimale d'une courbe ferm\'{e}e essentielle non p\'{e}riph\'{e}rique dans $\hat{Y}$. Il montre le th\'{e}or\`{e}me suivant~:
 
\begin{theorem} [Rafi]
Il existe une constante $C$ qui d\'{e}pend uniquement de $\delta$ et de la topologie de $X$ telle pour toute courbe $\alpha$ essentielle contenue dans $\hat{Y}$, on ait
\[ \
{C}^{-1} \, \lambda(\hat{Y}) \, \ell_{\mathrm{hyp}}(\alpha) \leq \ell_{\mathrm{flat}}(\alpha) \leq C \, \lambda(\hat{Y}) \, \ell_{\mathrm{hyp}} (\alpha).
\]
De plus,
\[
\frac{1}{2}\,\lambda(\hat{Y}) \leq \mathrm{Diam}(\hat{Y}) \leq C \lambda(\hat{Y})
\]
o\`{u} 
$\mathrm{Diam}(\hat{Y})$ est le diam\`{e}tre de $\hat{Y}$.
\end{theorem} 
Ce th\'{e}or\`{e}me signifie que la taille est le bon invariant num\'{e}rique pour comparer  longueur plate et  longueur hyperbolique.
Ce r\'{e}sultat de Rafi est le pas essentiel pour obtenir l'\'{e}nonc\'{e} sur la compactification. Celui-ci affirme que lorsque une suite de surfaces de translation converge vers une surface de Riemann stable, sur une partie \'{e}paisse, quitte \`{a} prendre une sous-suite, les 1--formes convergent vers une 1--forme non nulle \`{a} condition de renormaliser par la taille \footnote{Si on ne renormalise pas, on peut tr\`{e}s facilement imaginer une partie de la surface stable o\`{u} la forme limite est nulle.}. Comme le sugg\`{e}re la figure  
\ref{fig:2tores}
l'exemple classique consiste \`{a} coller ensemble deux tores plats d'aire respective 1 et $\varepsilon$ avec $\varepsilon$ tendant vers 0. 
Sur la partie droite de la figure  \ref{fig:2tores} la forme limite est nulle si on ne renormalise pas.  Par contre la m\'{e}trique hyperbolique est bien sur invariante par homoth\'{e}tie.

\par \smallskip
Le th\'{e}or\`{e}me suivant exprime pr\'{e}cis\'{e}ment ce comportement:

\begin{theorem} \label{thm:compactification}
Soit $(X_\tau, \omega_\tau)$ une suite de surfaces de translation dans une strate $\HHH(k_1, \dots, k_r)$ telle que $(X_\tau)$ converge vers une surface de Riemann stable $X_\infty$.  Consid\'{e}rons une composante irr\'{e}ductible $Y_{\infty, j}$ de $X_\infty$, soit $Y_{\tau, j} $ la composante \'{e}paisse correspondante dans $X_\tau$ et $\lambda (  Y_{\tau, j} )$ la taille de son repr\'{e}sentant plat. Notons
$\omega_{\tau, j} ={\lambda (  Y_{\tau, j} )}^{-1}\, {\omega_{\tau}}$. 
Alors une sous suite  de $(\omega_{\tau, j})_{\tau}$ restreinte \`{a} $\hat{Y}_{\tau, j}$ converge vers une 1--forme m\'{e}romorphe non nulle $\omega_{j}$ sur (le repr\'{e}sentant plat de) $Y_{\infty, j}$. 
Les z\'{e}ros et p\^{o}les de $\omega_{j}$ sont les limites de ceux de $\omega_{\tau', j}$ et eventuellement les noeuds.
\end{theorem}

La d\'{e}monstration de ce th\'{e}or\`{e}me n'est pas difficile une fois admis le th\'{e}or\`{e}me de Rafi. A quelques d\'etails pr\`es, il faut construire une triangulation  de $Y_{\tau, j}$ par liens de selles de longueur born\'{e}e apr\`{e}s avoir renormalis\'{e} par la taille. Ceci s'obtient par r\'{e}currence vu que le diam\`{e}tre de $Y_{\tau, j}$ est contr\^{o}l\'{e} par la taille. On a alors un ensemble compact de structures plates (vu que les longueurs des ar\^{e}tes de la triangulation sont born\'{e}es), il suffit alors de prendre une limite de ces structures.

\subsection{D\'{e}terminant du Laplacien pr\`{e}s du bord et th\'{e}or\`{e}me de Lundelius}

Comme nous l'avons d\'{e}j\`{a} expliqu\'{e}, une id\'{e}e fondamentale est d'\'{e}tablir un lien entre le d\'{e}terminant relatif du laplacien pour la m\'{e}trique plate et pour la m\'{e}trique hyperbolique.  Eskin, Kontsevich et Zorich donnent une estim\'{e}e sur la diff\'{e}rence entre ces deux quantit\'{e}s et montrent que leurs valeurs au bord de l'espace des modules sont \'{e}quivalentes, m\^{e}me si elles sont toutes deux non born\'{e}es. Dans toute la suite, nous notons $\ell_{\mathrm{plat}}(S)$ la longueur du plus petit lien de selles sur une surface de translation $S$. Le r\'{e}sultat est le suivant:
\begin{theorem} \label{thm:terme-erreur}
Soit $(X, \omega)$ une surface de translation et $X_0$ un point base dans la m\^{e}me strate que $X$. Alors
\[
\big{\vert} \log \det \Delta_{\textrm{plat}}(X, X_0) -  \log \det \Delta_{\textrm{hyp}}(X, X_0) \big{\vert}  = O \bigl(\vert \log \ell_{\textrm{plat}}(X)\vert \bigr).
\]
\end{theorem}
\par \smallskip
{\it Dans l'article original, les auteurs produisent des estim\'{e}es explicites. Les constantes sont universelles et ne d\'{e}pendent que du genre, du nombre de z\'{e}ros et de la surface de r\'{e}f\'{e}rence $X_0$.}
\par \medskip
 
Un th\'{e}or\`{e}me de Lundelius \cite{Lu} donne un d\'{e}veloppement asymptotique de $\log \det \Delta_{\textrm{hyp}}$ lorsque un certain nombre de courbes d\'{e}g\'{e}n\`{e}rent. 

\begin{theorem} (Lundelius)
Soit $(C_{\tau})$ une famille de surfaces hyperboliques de volume fini, de type topologique fix\'{e}, qui tend vers une courbe stable $C_{\infty}$, et $C_0$ une surface de r\'{e}f\'{e}rence de m\^{e}me type topologique que $C_\tau$. Alors
\[
-\log \vert \det \Delta_{\textrm{hyp}}(C_\tau, C_0) \vert = \sum_k \frac{\pi^2}{3 \ell_{\tau, k}} + O \bigl(-\log \ell_{\mathrm{hyp}}(C_\tau) \bigr) + O(1)
\]
lorsque $\tau$ tend vers l'infini, o\`{u} $ \ell_{\tau, k}$ sont les longueurs des g\'{e}od\'{e}siques hyperboliques pinc\'{e}es et 
$\ell_{\mathrm{hyp}}(C_\tau)$ est la longueur hyperbolique de la plus petite g\'{e}od\'{e}sique sur $C_\tau$.

\end{theorem}

Ce th\'{e}or\`{e}me est une g\'{e}n\'{e}ralisation de r\'{e}sultats de Wolpert \cite{Wo} et d'Osgood, Phillips et Sarnak \cite{OPS} dans le cas compact. 
Une \'{e}tape de la preuve du th\'{e}or\`{e}me \ref{cacasselabaraque} consiste \`{a} combiner le th\'{e}or\`{e}me \ref{thm:terme-erreur} avec le th\'{e}or\`{e}me de Lundelius pour obtenir un d\'{e}veloppement asymptotique de la fonction z\^{e}ta de $\Delta_{\mathrm{plat}}$. Les auteurs prouvent le th\'{e}or\`{e}me suivant~:

\begin{theorem} \label{thm:cylindres-plats}
Pour toute strate $\HHH(k_1, \dots, k_r)$, il existe une constante $M$ telle que 
\[-\log \vert \det \Delta_{\mathrm{plat}}(X, X_0) \vert = \frac{\pi}{3}\sum_{\substack {\mathrm{cylindres} \ \CCC \\ \mathrm{avec } \ h_{\CCC}/w_{\CCC} \geq M}} \frac{h_{\CCC}}{w_{\CCC}}  + O \bigl(-\log \ell_{\mathrm{plat}}(S) \bigr)
\]
o\`{u} $h_{\CCC}$ est la hauteur du cylindre $\CCC$, $w_{\CCC}$ sa circonf\'{e}rence et $X_0$ une surface de r\'{e}f\'{e}rence dans la strate $\HHH(k_1, \dots, k_r)$.
\end{theorem}

Les d\'{e}monstrations de ces th\'{e}or\`{e}mes sont tr\`{e}s techniques, calculatoires et d\'{e}velopp\'{e}es sur plusieurs dizaines de pages, nous n'allons \'{e}videmment pas les d\'{e}tailler. Nous donnons simplement quelques rep\`{e}res qui aideront le lecteur \`{a} comprendre la structure de la preuve. 
\par \medskip
La premi\`{e}re id\'{e}e est la formule de Polyakov \eqref{eq:Polyakov} qui permet d'exprimer la diff\'{e}rence 
\[
\log \det \Delta_{\mathrm{plat}}(S, S_0) - \log \det \Delta_{\mathrm{hyp}}(S,S_0)
\] 
sous forme int\'{e}grale. Nous avons expliqu\'{e} plus haut que l'on consid\`{e}re des d\'{e}terminants relatifs en r\'{e}gularisant les m\'{e}triques (plates et hyperboliques) aux voisinage des z\'{e}ros. Il faut donc comprendre la contribution de ces points. C'est l'objet de la partie 6 de l'article qui est essentiellement calculatoire.
Ensuite, il faut tirer profit du th\'{e}or\`{e}me  \ref{thm:compactification} et estimer \[
\log \det \Delta_{\mathrm{plat}}(S, S_0) - \log \det \Delta_{\mathrm{hyp}}(S,S_0)
\] 
lorsque les surfaces d\'{e}g\'{e}n\`{e}rent. On obtient alors une version forte du th\'{e}or\`{e}me \ref{thm:terme-erreur}.
Pour passer du th\'{e}or\`{e}me de Lundelius au th\'{e}or\`{e}me  \ref{thm:cylindres-plats},  il faut montrer que la contribution des anneaux euclidiens est n\'{e}gligeable et que 
$\vert \log \ell_{\mathrm{hyp}} \vert = O \bigl(\vert \log \ell_{\mathrm{plat}} \vert \bigr)$.
Le lemme cl\'{e} est ici la relation entre module et longueur de g\'{e}od\'{e}sique~:

\begin{lem} \label{lem:module}
Si $A$ un cylindre de module $M$, il existe une unique g\'{e}od\'{e}sique pour {\bf sa}  m\'{e}trique hyperbolique dont la longueur hyperbolique est $\pi/M$.
\end{lem}

Ceci est totalement classique\footnote{On pourra voir le livre de Hubbard \cite{Hu} o\`{u} les calculs de g\'{e}om\'{e}trie hyperbolique sont fort bien pr\'{e}sent\'{e}s.}. Lorsqu'un cylindre de grand module est contenu dans une surface hyperbolique, on a approximativement le m\^{e}me r\'{e}sultat pour la m\'{e}trique hyperbolique de la surface.  La subtilit\'{e} est que la m\'{e}trique hyperbolique du cylindre n'est pas la m\^{e}me que celle de la surface ambiante. Les estim\'{e}es pr\'{e}cises sont dues \`{a} Masur, Minsky et Wolpert.
Pour les cylindres plats, la circonf\'{e}rence du cylindre est donc approximativement  la longueur de la g\'{e}od\'{e}sique hyperbolique correspondante. 
Par contre, pour les anneaux, la longueur de la g\'{e}od\'{e}sique hyperbolique est de l'ordre du logarithme de la longueur plate de la courbe correspondante. Ainsi, la contribution des anneaux est comprise dans le terme d'erreur en $O \bigl( \vert \log \ell_{\mathrm{plat}}(S) \vert \bigr)$.

\subsection{D\'{e}terminant du Laplacien et constantes de Siegel-Veech}

La derni\`{e}re partie importante de l'article (section 9) prouve le th\'{e}or\`{e}me suivant:

\begin{theorem} \label{thm:final}
Soit $\HHH(k_1, \dots, k_r)$ une strate de 1--formes holomorphes, $\nu_1$ la mesure normalis\'{e}e de Masur-Veech associ\'{e}e et 
$C_{\textrm{aire}}(\nu)$ la constante de Siegel-Veech correspondante. Alors
\[
\int_{\HHH_1(k_1, \dots, k_r)} \Delta_{\mathrm{Teich}} \log \Delta_{\mathrm{plat}}(X, X_0) \,d\nu_1(X) = -\frac{4}{3} \pi^2 \, C_{\mathrm{aire}}(\nu)
\]
o\`{u} $X_0$ est n'importe quelle surface de r\'{e}f\'{e}rence dans la strate. 
\end{theorem}

Ce r\'{e}sultat combin\'{e} avec  le th\'{e}or\`{e}me de Riemann-Roch analytique donne bien entendu le r\'{e}sultat principal. 
\par \medskip
Donnons tout d'abord une id\'{e}e heuristique du lien entre les deux membres de l'\'{e}galit\'{e} ci-dessus. En utilisant la formule de Green, le terme de gauche se r\'{e}\'{e}crit comme une int\'{e}grale sur un voisinage du bord de la strate. Vu le th\'{e}or\`{e}me \ref{thm:cylindres-plats},
le terme dominant est une int\'{e}grale sur les cylindres de grands modules. Les constantes de Siegel-Veech mesurent elles aussi la contribution de (certains) cylindres de grands modules. On compte donc intuitivement la m\^{e}me chose. 
Essayons d'expliquer cela un peu plus pr\'{e}cis\'{e}ment. 
\par \medskip
Rappelons tout d'abord la formule classique de Siegel-Veech \cite{Ve3} ainsi quelques faits d\^{u}s \`{a} Eskin-Masur \cite{EM} utilisant les constantes de Siegel-Veech.
\par \medskip
Soit $(X, \omega)$ une surface de translation, $f : \R^2 \to \R$ une fonction continue \`{a} support compact et $V(X)$ l'ensemble des vecteurs d'holonomie des liens de selles sur une surface de translation $X$
\footnote{$V(X)$ peut aussi \^{e}tre l'ensemble des vecteurs d'holonomie des \^{a}mes des cylindres sur $X$.}, c'est-\`{a}-dire les vecteurs $\displaystyle \int_{\gamma} \omega$ o\`{u} $\gamma$ parcourt l'ensemble des liens de selles sur $X$.
On d\'{e}finit la transform\'{e}e de Siegel-Veech de $\vphantom{\bigl(}f$ par la formule
\[
\hat{f}(X) = \sum_{v \in V(X)} f(v)~;
\] 
$\hat{f}$ est bien d\'{e}finie car $V(X)$ est un ensemble discret et $f$ est \`{a} support compact. 
\begin{lem} [Formule de Siegel-Veech]
Il existe une constante $\hat{c}\,(\nu)$ telle que, pour tout $f$ continue \`{a} support compact sur $\R^2$,
\[
\int_{\HHH_1(k_1, \dots , k_r)} \hat{f}(X) \, d\nu_1 = \hat{c}\,(\nu) \, \int_{\R^2} f(x,y) dx dy.
\]
\end{lem}

D'apr\`es \cite{EM}, la constante $\hat{c}\,(\nu)$ est exactement celle qui intervient dans l'asymptotique quadratique du nombre de liens de selles sur une surface g\'{e}n\'{e}rique. De plus, la formule est aussi valable pour des fonctions plus g\'{e}n\'{e}rales, comme par exemple les fonctions caract\'{e}ristiques de boules.
\par \medskip
La d\'{e}monstration de la formule est tr\`{e}s simple.
L'application \[ 
f \mapsto \int_{\HHH_1(k_1, \dots , k_r)} \hat{f}(X)\, d\nu_1
\]
est une forme lin\'{e}aire continue positive donc une mesure. Comme $V(X)$ est \'{e}quivariant sous l'action de $\textrm{SL}_2(\R)$, cette mesure est invariante par cette action. C'est donc une combinaison de la mesure de Dirac en $(0,0)$ et de la mesure de Lebesgue. On voit, en testant sur des fonctions ad hoc, que c'est un multiple de la mesure de Lebesgue, d'o\`{u} le r\'{e}sultat. 
\par \medskip
Voici un exemple d'application de cette formule:
\begin{lem} \label{lem:constantes-bord}
Il existe une constante $C >0$ telle que pour tout $\varepsilon >0$,
\[
\nu_1 \, \bigl{\{}(X,\omega) \,\,\textrm{telles que}\,\, \ell_{\mathrm{plat}} (X) < \varepsilon \bigr{\}} \leq C \varepsilon^2.
\]
\end{lem}

Nous donnons une id\'{e}e de la preuve car elle est simple et instructive. Notons $N(X,L)$ le nombre de liens de selles de longueur au plus $L$. On rappelle que cette quantit\'{e} croit quadratiquement quand $L$ tend vers l'infini, regardons son comportement quand $L$ tend vers 0. 
La formule de Siegel Veech appliqu\'{e}e \`{a} la fonction indicatrice de la boule de $\R^2$ centr\'{e}e en $0$ de rayon $\varepsilon$ nous dit que 
\[
\int_{\HHH_1(k_1, \dots, k_r)} N(X, \varepsilon) \, d\nu_1(X) = \hat{c}(\nu) \, \pi \varepsilon^2
\]
o\`{u} $\hat{c}(\nu)$ est la constante de Siegel-Veech du probl\`{e}me de comptage des liens de selles.
Vu que lorsque $\ell_{\mathrm{plat}} (X) < \varepsilon$ on a $N(X, \varepsilon) >1$, on obtient le r\'{e}sultat du lemme. 
\par \medskip
Une cons\'{e}quence de ce lemme est que pour tout $\beta <2$, $\ell_{\mathrm{plat}}^{-\beta}$ est int\'{e}grable. 
Eskin-Masur montrent \'{e}galement le r\'{e}sultat suivant:
\begin{theorem} [Eskin-Masur] \label{thm:EM}
Pour tout $\beta$ tel que $1 < \beta <2$ et pour tout $\varepsilon$ assez petit, on a
\[
N(X, \varepsilon)  = O \bigl(\ell_{\mathrm{plat}}(X)^{-\beta} \bigr).
\]
\end{theorem}

La preuve de ce r\'{e}sultat (\cite{EM} Th. 5.1) est nettement plus d\'{e}licate que celle du lemme \ref{lem:constantes-bord}.
\par \bigskip
Nous allons maintenant expliquer deux lemmes simples qui donnent un aper\c{c}u de la strat\'{e}gie de preuve du th\'{e}or\`{e}me \ref{thm:final}. Fixons un r\'{e}el $K$ qui sera toujours suppos\'{e} suffisamment grand, et notons $\mathrm{Cyl}_K(X)$  l'ensemble des cylindres de modules au moins $K$ sur $X$. 
Pour $K$ assez grand, les cylindres de  module au moins $K$ sont disjoints (car les g\'{e}od\'{e}siques hyperboliques associ\'{e}es sont simultan\'{e}ment courtes). 
Ainsi $\mathrm{Cyl}_K(X)$ contient au plus $3g-3 +r$ \'{e}l\'{e}ments. 
On d\'{e}finit
\[
\ell_K(X) = \min \bigl{\{} w(C),\, C \in \mathrm{Cyl}_K(X)^{\vphantom{A}}\bigl{\}}.
\]
On donnera 1 comme valeur \`{a} $\ell_K(X)$ lorsque l'ensemble est vide. 
Il est bien clair que \[
\ell_K(X) \geq \ell_{\mathrm{plat}} (X)
\] vu qu'un cylindre plat est bord\'{e} par des liens de selles. 
Soit $\chi_{\varepsilon}$ la fonction caract\'{e}ristique de l'ensemble des surfaces de translation $(X, \omega)$ tels que $\ell_K(X) \geq \varepsilon$. On d\'{e}finit en premier lieu une fonction auxiliaire qui sert \`{a} r\'{e}gulariser. Soit $\eta  \colon \textrm{SL}_2(\R) \to \R$ une fonction r\'{e}guli\`{e}re positive $\textrm{SO}_2(\R)$ invariante telle que 
$\displaystyle \int_{\textrm{SL}_2(\R)} \eta(g) dg = 1$, de support contenu dans la couronne des \'{e}l\'{e}ments $g$ tels que $1/2 < \vert \vert g \vert \vert < 2$. La norme sur $\vphantom{\bigl)}\textrm{SL}_2(\R)$ est la norme induite par la norme euclidienne de $\R^2$. On d\'{e}finit
\[
f_{\varepsilon} (X) = \int_{\textrm{SL}_2(\R)} \eta(g)\chi_\varepsilon(gX) \, dg.
\]
La fonction $f_{\varepsilon}$ satisfait les propri\'{e}t\'{e}s suivantes:
\begin{enumerate}
\item $f_{\varepsilon} (X) = 0 \mbox{ si } \ell_K(X) \leq \varepsilon/2$
\item $f_{\varepsilon} (X) = 1 \mbox{ si } \ell_K(X) \geq 2\varepsilon$
\item $f_{\varepsilon}$ est r\'{e}guli\`{e}re le long des disques de Teichm\"{u}ller et son gradient ainsi que son laplacien sont born\'{e}s ind\'{e}pendamment de $\varepsilon$.
\end{enumerate}

Autrement dit, le gradient de  $f_{\varepsilon}$ est concentr\'{e} sur une couronne 
\[
\AAA_{\varepsilon} = \{X \mbox{ tel que } \varepsilon/2 \leq \ell_K(X) \leq 2\varepsilon\}
\]
qui est bien s\^{u}r un voisinage du bord et 
$f_\varepsilon (X)$ tend vers 1 lorsque $\varepsilon$ tend vers 0.
\par \medskip
Posons
\[
\psi^K(X) = \sum_{C \in \mathrm{Cyl}_K(X)} \bigl(\mathrm{Mod}(C) - K\bigr).
\]
Par construction, la fonction $\psi^K$ est une fonction continue et r\'{e}guli\`{e}re par morceaux. 
Pour se rapprocher un peu du probl\`{e}me de comptage intervenant dans la d\'{e}finition des constantes de Siegel-Veech, il est important de ne consid\'{e}rer que les cylindres qui sont parall\`{e}les \`{a} celui dont l'\^{a}me est la plus courte. 
On note $\widetilde{\mathrm{Cyl}}_K(X)$ l'ensemble de ces cylindres et on pose
\[
\widetilde{\psi}^K(X) = \sum_{C \in \widetilde{\mathrm{Cyl}}_K(X)} \bigl(\mathrm{Mod}(C) - K\bigr).
\]
Les fonctions $\psi^K$ et $\widetilde{\psi}^K$ ont un comportement analogue pr\`{e}s du bord. Ce point est technique mais important: c'est exactement \`{a} ce niveau que l'hypoth\`{e}se de r\'{e}gularit\'{e} de la mesure intervient dans la preuve g\'{e}n\'{e}rale. Pour les strates, Eskin-Masur-Zorich montrent que l'ensemble des surfaces qui ont deux liens de selles non homologues de longueur au plus $\varepsilon$ est de mesure $\varepsilon^4$. Le lemme \ref{lem:constantes-bord} dit que l'ensemble des surfaces avec un lien de selles de longueur au plus $\varepsilon$ a une mesure au plus de l'ordre de $\varepsilon^2$. Ainsi des surfaces g\'{e}n\'{e}riques proches du bord ont uniquement des liens de selles courts parall\`{e}les. On montre alors: 
\begin{lem} \label{premiersimple} Avec les notations pr\'{e}c\'{e}dentes, on a
\[
\int_{\HHH_1(k_1, \dots, k_r)}  \Delta_{\mathrm{Teich}}\log\Delta_{\mathrm{plat}}(X,X_0)\, d\nu_1 = \frac{\pi}{3}\, \lim_{\varepsilon \to  0}\, \int_{\HHH_1(k_1, \dots, k_r)} \nabla_{\mathrm{Teich}} \widetilde{\psi}^K(X) .\nabla_{\mathrm{Teich}} f_{\varepsilon}\, d\nu_1.
\]
\end{lem}

Donnons un sch\'{e}ma preuve pour comprendre comment interviennent les arguments introduits plus haut.
Posons $f = \log\Delta_{\mathrm{plat}}(X,X_0)$ pour simplifier les notations. 
Par la formule de Green, on a 
\[
\int_{\HHH_1(k_1, \dots, k_r)}  \Delta_{\mathrm{Teich}}f \, d\nu_1 = \lim_{\varepsilon \to 0} \int_{\HHH_1(k_1, \dots, k_r)} f_\varepsilon \,  \Delta_{\mathrm{Teich}}f \, d\nu_1 = 
\lim_{\varepsilon \to 0} \,  \int_{\HHH_1(k_1, \dots, k_r)} f  \, \Delta_{\mathrm{Teich}}f _\varepsilon \, d\nu_1.
\]
On utilise ici le th\'{e}or\`{e}me \ref{thm:cylindres-plats} pour estimer $f$ pour obtenir le d\'{e}veloppement
\[
f(X) = -\frac{\pi}{3}\psi^K(X) + O \bigl(\log (\ell_{\mathrm{plat}} (X) \bigr).
\]
 
Rappelons que $\Delta_{\mathrm{Teich}}f_{\varepsilon}$ est born\'{e} et que son support est contenu dans la couronne $\AAA_{\varepsilon}$. D'apr\`{e}s le lemme \ref{lem:constantes-bord}
 et le fait que $\ell_K(X) \geq \ell_{\mathrm{plat}} (X)$, on a $\nu_1(\AAA_{\varepsilon}) = O(\varepsilon^2).$ De plus, le  lemme pr\'{e}-cit\'{e} entraine que la fonction $X \rightarrow \log\bigl(\ell_{\mathrm{plat}} (X) \bigr)$ est int\'{e}grable sur $\HHH_1(k_1, \dots, k_r)$. Ainsi, en appliquant le th\'{e}or\`{e}me de convergence domin\'{e}e, on obtient
 \[
 \lim_{\varepsilon \to 0}\,  \int_{\HHH_1(k_1, \dots, k_r)} \vert \log \ell_{\mathrm{plat}} \vert \, \Delta_{\mathrm{Teich}} f_{\varepsilon} = 0.
 \]
 
 On obtient donc
 \[
 \int_{\HHH_1(k_1, \dots, k_r)}  \Delta_{\mathrm{Teich}}\log\Delta_{\mathrm{plat}}(X,X_0) \, d\nu_1 = \frac{\pi}{3} \, \lim_{\varepsilon \to  0}\,  \int_{\HHH_1(k_1, \dots, k_r)} \nabla_{\mathrm{Teich}} \psi^K(X) .\nabla_{\mathrm{Teich}} f_{\varepsilon} \, d\nu_1.
 \]

Reste \`{a} comparer $\psi^K$ et $\widetilde{\psi}^K$ au voisinage du bord ce qui se fait sans mal avec le discours pr\'{e}c\'{e}dent.
\par \medskip
Nous passons maintenant au second lemme qui exprime la constante de Siegel-Veech sous forme int\'{e}grale. 
Notons $\mathrm{Cyl}(X, \varepsilon)$ les cylindres plats sur $X$ qui ont une \^{a}me de longueur plate entre $\varepsilon/2$ et $\varepsilon$. 

\begin{lem} \label{secondsimple}
Si on pose
\[
\widetilde{N}^K_{\mathrm{aire}}(X, \varepsilon) = \sum_{C \in \widetilde{\mathrm{Cyl}}_K (X) \,\cap \, \mathrm{Cyl}(X, \varepsilon)} \mathrm{Aire}(C)
\]
on a 
\[
C_{\mathrm{aire}}(\nu) = \lim_{\varepsilon \to 0} \, \frac{4}{3 \pi \varepsilon^2}\,  \int_{\HHH_1(k_1, \dots, k_r)} \widetilde{N}^K_{\mathrm{aire}}(X, \varepsilon) \, d\nu_1(X).
\]
\end{lem}

Montrons le r\'{e}sultat en rempla\c{c}ant $\widetilde{\mathrm{Cyl}}_K $ par $\mathrm{Cyl}_K$ ce qui ne change presque rien.  
Notons $\hat{N}_{\mathrm{aire}}(X, \varepsilon)$  la contribution d'une couronne, c'est-\`{a}-dire $ N_{\mathrm{aire}}(X, \varepsilon) - N_{\mathrm{aire}}(X, \varepsilon/2)$. 
Une fois de plus, la formule de Siegel-Veech nous permet d'\'{e}crire $C_{\mathrm{aire}}(\nu)$ sous forme int\'{e}grale~: pour tout $\varepsilon >0$,
\begin{equation} \label{eq:SV}
C_{\mathrm{aire}}(\nu) = \frac{4}{3\pi \varepsilon^2} \int_{\HHH_1(k_1, \dots, k_r)} \hat{N}_{\mathrm{aire}}(X, \varepsilon)\, d\nu_1(X).
\end{equation}

Si $C$ est un cylindre dans $\mathrm{Cyl}(X,\varepsilon) \setminus \mathrm{Cyl}_K(X)$, 
\[
\mathrm{Aire(C)}= \dfrac{h}{w} \times w^2 < K\varepsilon^2.
\]
Par le th\'{e}or\`{e}me \ref{thm:EM}, $N(X, \varepsilon)  = O\bigl(\ell_{\mathrm{plat}}(X)^{-\beta} \bigr)$
pour $1 < \beta <2$ donc 
\[
\hat{N}_{\mathrm{aire}}(X,\varepsilon) - N^K_{\mathrm{aire}}(X, \varepsilon) \leq \kappa \, K\varepsilon^2\,  \ell_{\mathrm{plat}}(X)^{-\beta}
\]
 o\`{u} $\kappa$ est une constante ind\'{e}pendante de $\varepsilon$.
On int\`{e}gre cette in\'{e}galit\'{e} sur l'ensemble des surfaces de translation $(X,\omega)$ telles que  $\ell(X) < \varepsilon$. Pour le membre de gauche cela revient \`{a} int\'{e}grer sur la strate tout enti\`{e}re. Ainsi
\[
\frac{1}{\varepsilon^2}\int_{\HHH_1(k_1, \dots, k_r)} \bigl( \hat{N}_{\mathrm{aire}}(X,\varepsilon) - N^K_{\mathrm{aire}}(X, \varepsilon)\bigr)\, d\nu_1(X) \leq \kappa K  \int_{\{ (X,\omega) \, \vert \, \ell_{\mathrm{plat}} (X) < \varepsilon    \}} \ell_{\mathrm{plat}}(X)^{-\beta} \, d\nu_1(X).
\]
L'int\'{e}grabilit\'{e} de $\ell_{\mathrm{plat}}^{-\beta} $ et le lemme \ref{lem:constantes-bord} impliquent que 
\[
\lim_{\varepsilon \to 0} \, \frac{1}{\varepsilon^2} \, \int_{\HHH_1(k_1, \dots, k_r)} \bigl( \hat{N}_{\mathrm{aire}}(X,\varepsilon) - N^K_{\mathrm{aire}}(X, \varepsilon) \bigr) \, d\nu_1(X) =0.
\]
Dans l'\'{e}quation \eqref{eq:SV}, on peut donc remplacer $\hat{N}_{\mathrm{aire}}(X, \varepsilon)$ par $N^K_{\mathrm{aire}}(X, \varepsilon)$ ce qui prouve le lemme.
 \par \medskip
Pour conclure, on voit bien que les lemmes \ref{premiersimple} et \ref{secondsimple} permettent de rapprocher la moyenne de $\Delta_{\mathrm{Teich}} \log\det\Delta_{\mathrm{flat}}$ et la constante de Siegel-Veech $C_{\mathrm{aire}}$. Quelques efforts suppl\'{e}mentaires et pas mal de technique sont n\'{e}cessaires pour obtenir une preuve compl\`{e}te du Th\'{e}or\`{e}me \ref{thm:final}. 

\bigskip

{\bf Remerciements.} Nous remercions tout d'abord les auteurs qui ont r\'epondu \`a nos questions avec beaucoup de patience et de gentillesse. Ce texte n'aurait pas vu le jour sans C. Boissy et E. Lanneau, merci \`a eux deux. Nous remercions enfin G. Merlet, J. H. Hubbard, E. Russ, C.~Soul\'e pour des discussions fort utiles.

\end{document}